\documentclass[ejs]{imsart}

\usepackage{verbatim}
\usepackage{amsmath,amsthm,amssymb}
\usepackage{graphicx}
\usepackage{amsthm,amsmath}

\usepackage[usenames]{color}
\RequirePackage[colorlinks,citecolor=blue,urlcolor=blue]{hyperref}

\usepackage{algorithm,algorithmic}

\startlocaldefs

\def\deq{\triangleq}

\def\R{{\mathbb R}}
\def\N{{\mathbb N}}
\def\E{{\mathbb E}}
\def\var{\operatorname{Var}}
\def\Pr{{\mathbb P}}

\def\cF{{\cal F}}
\def\cL{{\cal L}}

\def\chr{\varphi}
\def\bin{{\mathbb B}}

\def\wh#1{\widehat{#1}}

\def\ave#1{\langle #1 \rangle}

\def\Sum{\displaystyle\sum}

\newtheorem{theorem}{Theorem}
\newtheorem{lemma}{Lemma}
\newtheorem{remark}{Remark}

\endlocaldefs

\begin{document}

\begin{frontmatter}

% "Title of the paper"
\title{A recursive procedure for density estimation on the binary hypercube\protect\thanksref{T1}}
\runtitle{Recursive density estimation on the binary hypercube}

\begin{aug}
\author{Maxim Raginsky\ead[label=e1]{maxim@illinois.edu}}

\address{Department of Electrical and Computer Engineering\\
and Coordinated Science Laboratory\\
University of Illinois\\
Urbana, IL 61801\\
\printead{e1}}

\author{Jorge Silva\ead[label=e2]{jorge.gomes.da.silva@gmail.com}}

\address{Department of Electrical and Computer Engineering\\
Duke University\\
Durham, NC 27708\\
\printead{e2}}

\author{Svetlana Lazebnik\ead[label=e4]{slazebni@illinois.edu}}

\address{Department of Computer Science\\
University of Illinois\\
Urbana, IL 61801\\
\printead{e4}}

\author{Rebecca Willett\ead[label=e3]{willett@duke.edu}}

\address{Department of Electrical and Computer Engineering\\
Duke University\\
Durham, NC 27708\\
\printead{e3}}

\thankstext{T1}{This work was supported by NSF CAREER Award CCF-06-43947, DARPA Grant HR0011-07-1-003, and ARO Grant W911NF-09-1-0262.}

\end{aug}

\runauthor{Raginsky et al.}

\begin{abstract}
	This paper describes a recursive estimation procedure for multivariate binary densities (probability distributions of vectors of Bernoulli random variables) using orthogonal expansions. For $d$ covariates, there are $2^d$ basis coefficients to estimate, which renders conventional approaches computationally prohibitive when $d$ is large. However, for a wide class of densities that satisfy a certain sparsity condition, our estimator runs in probabilistic polynomial time and adapts to the unknown sparsity of the underlying density in two key ways: (1) it attains near-minimax mean-squared error for moderate sample sizes, and (2) the computational complexity is lower for sparser densities. Our method also allows for flexible control of the trade-off between mean-squared error and computational complexity.
\end{abstract}



\tableofcontents

\end{frontmatter}

\section{Introduction}
\label{sec:intro}

This paper considers the problem of estimating a multivariate binary density from a number of independent observations. That is, we have $n$ observations of the form $X_i \in \{0,1\}^d$ which are independent and identically distributed (i.i.d.) samples from a population with a probability density $f$ (with respect to the counting measure on the $d$-dimensional {\em binary hypercube} $\{0,1\}^d$). We wish to estimate $f$ on the basis of these observations. Multivariate binary data arise in a variety of applications:
\begin{itemize}
\item {\em Biostatistics.} Each $X_i$ could represent a biochemical profile of a bacterial strain, where every component is a ``yes-no" indicator of a presence of a particular biochemical marker \cite{bacterialTaxonomy,communityDNA}. A recent paper \cite{ShmZha02} proposed a methodology for representing gene expression data using binary vectors. More classical scenarios include recording the occurrence of a given symptom or a medical condition in a patient over time \cite{dental} or outcomes of a series of medical tests \cite{AitAit76}.
\item {\em Quantitative methods in social sciences.} Each $X_i$ could
  represent a respondent in a survey or a panel, where every component
  is a ``yes-no" answer to a question \cite{Car07}, describe a voting record of
  a legislator, or correspond to co-occurrences of events in social networks \cite{silvaHypergraph}.
\item {\em Artificial intelligence.} Each $X_i$ could represent a user query to a search engine or a database, where every component corresponds to the presence or absence of a particular keyword \cite{GhaHel06}, or an image stored on a website like Flickr\footnote{{\tt http://www.flickr.com}}, where every component corresponds to a user-supplied tag from a given list.
\end{itemize}
Many situations involving multivariate binary data have the following features: (1) the number of covariates (or the dimension of the hypercube) $d$ is such that the number of possible values each observation could take ($2^d$) is much larger than the sample size $n$; (2) there is a ``clustering effect" in the population, meaning that the shape of the underlying density is strongly influenced mainly by a small number of constellations of the $d$ covariates. For example, a particular class of bacterial strains may be reliably identified by looking at a particular subset of the biomarkers; there may be several such classes in the population of interest, each associated with a distinct subset of biomarkers. Similarly, when working with panel data, it may be the case that the answers to some specific subset of questions are highly correlated among a particular group of the panel participants, and the responses of these participants to other questions are nearly random; moreover, there may be several such distinct groups in the panel. 

These considerations call for a density estimation procedure that can effectively cope with ``thin" samples (i.e.,~those samples for which $n < 2^d$) in terms of both estimation error and computational complexity, and at the same time automatically adapt to the possible clustering in the population, in the sense described above. We take the minimax point of view, where we assume that the unknown density $f$ comes from a particular function class $\cF$ and seek an estimator that exactly or approximately attains the minimax mean-squared error
$$
R^*_n(\cF) = \inf_{\widehat{f}} \sup_{f \in \cF} \E \| \widehat{f} - f \|^2_{L^2},
$$
where the infimum is over all estimators based on $n$ i.i.d.\ samples from $f$. We will choose the class $\cF$ to model the ``constellation effects" via a certain {\em sparsity} condition. Our choice of the $L^2$ risk, as opposed to other measures of risk more commonly used in density estimation, such as Hellinger, Kullback--Leibler or total variation risks, is dictated by the fact that the sparsity condition mentioned above is most naturally stated in a Hilbert space framework, which in turn facilitates the design of our estimation procedure, as well as the derivation of both upper and lower bounds on $R^*_n(f)$. We refer the reader to several other works on density estimation that use $L^2$ risk \cite{LiaKri85,CheKriLia89,DJKP96,YangBarron,HPKP97,HalKerPic98,Efromovich}. We also note that, because the Euclidean $L^1$ norm (which for probability densities gives the total variation risk) dominates the Euclidean $L^2$ norm, and because the square of the Kullback--Leibler divergence dominates the total variation distance \cite{CoverThomas}, lower bounds on the squared $L^2$ risk automatically translate into lower bounds on the squared $L^1$ (total variation) risk and on the Kullback--Leibler risk.

Because of the host of applications in which multivariate
binary data naturally arise, several authors have investigated
algorithms for estimation of their probability densities
(see, e.g.,~\cite{AitAit76,OttKro76,LiaKri85,CheKriLia89}). However, existing
approaches either have very slow rates of error convergence
or are computationally prohibitive when the number of covariates is
very large. For example, the kernel density estimation scheme proposed by Aitchison and Aitken \cite{AitAit76} has computational complexity $O(nd)$ (where $n$ is the number of observations and $d$
the number of covariates), yet its squared $L^2$ error
decays at the rate $O(n^{-4/(4+d)})$ \cite{Sim95}, which is disastrously slow for
large $d$. In contrast, orthogonal series methods, which can
potentially achieve near-minimax error decay rates, require the
estimation of $2^d$ basis coefficients and do not easily admit
computationally tractable estimation methods for very large
$d$. For instance, using the Fast Walsh--Hadamard Transform to estimate the coefficients of a density in the Walsh basis (see below) using $n$ samples requires $O(nd2^d)$ operations (see Appendix~B in \cite{OttKro76} and references therein).

In this paper we present a computationally tractable orthogonal series estimation method based on recursive block thresholding of empirical Walsh coefficients. In particular, the proposed method entails
recursively examining empirical estimates of whole {\em blocks} of the $2^d$ different Walsh coefficients. At each stage of the algorithm, the overall weight of basis coefficients computed at previous stages is used to decide
which remaining coefficients are most likely to be significant or
insignificant, and computing resources are allocated accordingly. It is
shown that this decision is accurate with high probability, so that insignificant coefficients are not estimated, while the
significant coefficients are. This approach is similar in spirit to the algorithm of Goldreich and Levin \cite{GolLev89}, originally developed for applications to cryptography and later adopted by Kushilevitz and Mansour \cite{KusMan93,Man94} to the problem of learning Boolean functions using membership queries. Although there are significant differences between the problems of density estimation and function learning which are reflected in our estimation procedure, our algorithm inherits the computational tractability of the Goldreich--Levin scheme: in particular, it runs in probabilistic polynomial time.

\sloppypar The proposed estimator adapts to unknown sparsity of the underlying density in two distinct ways. First, it is near-minimax optimal for ``moderate'' sample sizes $d \preceq n \preceq 2^{2d/p}$, with an $L^2$ error
decay rate of $O(2^{-d}(d/n)^{2r/(2r+1)})$, where $p \in (0,1]$ is a measure of sparsity and $r=1/2-1/p$. Moreover, the computational complexity of our algorithm is automatically lower for sparser densities. Sparsity has been recently recognized as a crucial
enabler of accurate estimation in ``big-$d$, small-$n$'' type problems \cite{Can06,CanTao06}. Specifically for densities on the binary hypercube, sparsity in the Walsh basis has a natural qualitative interpretation that the shape of the density is influenced mainly by a small number of constellations of the covariates. For example, if the components of a multidimensional binary vector represent positive/negative outcomes in a series of medical tests, it is often the case that the outcomes of certain small constellations of tests play the determining role in the diagnosis.

There are several different series expansions on the binary
hypercube presented in the literature, including the Rademacher--Walsh
orthogonal series (see Appendix~A in \cite{OttKro76} and references therein) and the Bahadur expansion
\cite{Bah61,bacterialTaxonomy}.  We focus in this paper on the Walsh system,
which is derived from Fourier analysis on finite groups (see, e.g., Chap.~4 of Tao and Vu \cite{TaoVu06}), for two 
reasons. First, the coefficients of a particular function in the Walsh system give us information about the influence of the various subsets of the $d$ variables on the value of the function \cite{Tal94a,DFKO07}. Second, the Walsh functions of a length-$d$ vector
can be factorized into products of Walsh functions of
multiple shorter vectors with lengths summing to $d$; this is detailed
in Section~\ref{ssec:walsh}. This factorization is central to the
efficiency of the proposed coefficient estimation method. The Walsh system is widely used in the context of learning Boolean functions \cite{Man94}, as well as in harmonic analysis of real-valued functions on the binary hypercube \cite{Tal94a,DFKO07}.

\subsection{Organization of the paper}
\label{ssec:paper}

The remainder of the paper is organized as follows. Section~\ref{sec:prelims} contains the preliminaries on notation, the Walsh system, and sparsity classes on the binary hypercube. Next, in Section~\ref{sec:recursive} we describe the motivation behind the thresholding methods in orthogonal series estimation on the binary hypercube, introduce our recursive thresholded estimator, and analyze its MSE and computational complexity. The theorems of Section~\ref{sec:recursive} are proved in Section~\ref{sec:proofs}.  Some illustrative simulation results are given in Section~\ref{sec:simulation}. The contributions of the paper are summarized in Section~\ref{sec:conclusion}. Finally, some technical results are relegated to the appendices.

\section{Preliminaries}
\label{sec:prelims}

\subsection{Notation}
\label{ssec:notation}

The basic set $\{0,1\}$ will be denoted by $\bin$. For any integer $k > 1$, the components of binary strings $x \in \bin^k$ will be denoted by $x^{(j)}$, $1 \le j \le k$: for any $x \in \bin^k$, we have $x = (x^{(1)},\ldots,x^{(k)})$. We will use juxtaposition to denote concatenation of strings: if $y \in \bin^k$ and $z \in \bin^l$, then $yz \in \bin^{k+l}$ is the string $x = (y^{(1)},\ldots,y^{(k)},z^{(1)},\ldots,z^{(l)})$. For any $0 \le k \le d$, we will define the {\em prefix} mapping $\pi_k : \bin^d \to \bin^k$ and the {\em suffix} mapping $\sigma_k : \bin^d \to \bin^{d-k}$ by
\begin{align}\label{eq:pref_suf}
\pi_k(x) \deq (x^{(1)},\ldots,x^{(k)}), \quad \sigma_k(x) \deq (x^{(k+1)},\ldots,x^{(d)}),
\end{align}
so that $x = \pi_k(x)\sigma_k(x)$ for any $x \in \bin^d$ (note that
both $\pi_0$ and $\sigma_d$ return an empty string). Whenever we deal
with vectors whose components are indexed by the elements of $\bin^k$
for some $k$, we will always assume that the components are arranged
according to the lexicographic ordering of the binary strings in
$\bin^k$. Given two real numbers $a,b$, we let $a \wedge b$ denote
$\min \{a,b\}$. Also, throughout the paper, $C$ is used to denote a
generic constant whose value may change from line to line; specific
absolute constants will be denoted by $C_1,C_2$, etc. 

Throughout the paper, we let $M \equiv 2^d$.

\subsection{The Walsh system}
\label{ssec:walsh}

For any integer $k \ge 1$, denote by $\mu_k$ the counting measure on $\bin^k$ and endow the space of functions $f : \bin^k \to \R$ with the structure of the real Hilbert space $L^2(\mu_k)$ via the standard inner product
$$
\ave{f,g} \deq \sum_{x \in \bin^k} f(x)g(x).
$$
The {\em Walsh system} (see references in the Introduction) in $L^2(\mu_k)$ is an orthonormal system $\Phi_k = \{\chr_s : s \in \bin^k\}$, defined by
\begin{equation}\label{eq:walsh}
\chr_s(x) \deq \frac{1}{2^{k/2}}(-1)^{s \cdot x}, \qquad \forall x \in \bin^k
\end{equation}
where $s \cdot x \deq \sum^k_{j=1} s^{(j)}x^{(j)}$. Hence, any $f \in L^2(\mu_k)$ has the Fourier--Walsh expansion
$$
f = \sum_{s \in \bin^k} \theta_s \chr_s,
$$
where $\theta_s \deq \ave{f,\chr_s}$, $s \in \bin^k$. To keep the notation simple, we will not explicitly mark the underlying dimension when working with the Walsh functions. When $k = d$, we will write $\Phi$ instead of $\Phi_d$.

For any $k$, the Walsh system $\Phi_k$ is a tensor product basis induced by $\Phi_1 = \{\chr_0,\chr_1\}$, where
$$
\chr_0(x) = \frac{1}{\sqrt{2}} \qquad \mbox{and} \qquad \chr_1(x) = \frac{1}{\sqrt{2}} (-1)^x
$$
for any $x \in \bin$. That is, for any $k \ge 1$, any $\chr_s \in \Phi_k$ has the form
$$
\chr_s = \chr_{s^{(1)}} \otimes \chr_{s^{(2)}} \otimes \ldots \otimes \chr_{s^{(k)}},
$$
which means that
$$
\chr_s(x) = \prod^k_{i=1} \chr_{s^{(i)}}(x^{(i)}), \qquad \forall x \in \bin^k.
$$
This generalizes to the following useful factorization property of the Walsh functions: for any $k \ge 1$ and any $l \le k$, we have
\begin{equation}\label{eq:walsh_factorize}
\chr_s = \chr_{\pi_{l,k}(s)} \otimes \chr_{\sigma_{l,k}(s)}, \qquad \forall s \in \bin^k
\end{equation}
where $\pi_{l,k}$ and $\sigma_{l,k}$ denote the prefix and the suffix mappings defined on $\bin^k$ analogously to \eqref{eq:pref_suf}. This means that, for products of functions on disjoint subsets of the $d$ variables, the Fourier--Walsh coefficients also have the product form.

\subsection{Sparsity and weak-$\ell^p$ balls}
\label{ssec:sparsity}

Our interest lies with functions whose Fourier--Walsh representations satisfy a certain sparsity constraint. Given a function $f$ on $\bin^d$, let $\theta(f)$ denote the vector of its Fourier--Walsh coefficients. We will assume that the components of $\theta(f)$ decay according to a power law. Formally, let $\theta_{(1)},\ldots,\theta_{(M)}$, where $M = 2^d$, be the components of $\theta(f)$ arranged in decreasing order of magnitude:
$$
|\theta_{(1)}| \ge |\theta_{(2)}| \ge \ldots \ge |\theta_{(M)}|.
$$
Given some $0 < p < \infty$, we say that $\theta(f)$ belongs to the Marcinkiewicz, or {\em weak-$\ell^p$}, ball of radius $R$ \cite{BerLof76,Joh94}, and write $\theta(f) \in w\ell^p(R)$, if
\begin{equation}\label{eq:power_law}
|\theta_{(m)}| \le R \cdot m^{-1/p}, \qquad  1 \le m \le M.
\end{equation}
It is not hard to show that the Fourier--Walsh coefficients of any
probability density on $\bin^d$ are bounded by $1/\sqrt{M}$. With this
in mind, let us define the function class
\begin{align}\label{eq:function_class}
	\cF_d(p) \deq \left\{ f : \bin^d \to \R: \theta(f) \in
  w\ell^p(1/\sqrt{M}) \right\}
\end{align}
We are particularly interested in the case $0 < p \le 1$.

We will need approximation properties of weak-$\ell^p$ balls as listed, e.g., in~\cite{CanTao06}. The basic fact is that the power-law condition (\ref{eq:power_law}) particularized to the elements of $\cF_d(p)$ is equivalent to the concentration estimate
\begin{equation}\label{eq:concentration}
\left| \left\{ s \in \bin^d : |\theta_s|^2 \ge \lambda \right\} \right| \le \left(\frac{1}{M\lambda}\right)^{p/2}
\end{equation}
valid for all $\lambda > 0$. Additionally, for any $1 \le k \le M$, let $\theta_k(f)$ denote the vector $\theta(f)$ with $\theta_{(k+1)},\ldots,\theta_{(M)}$ set to zero. Then it follows from (\ref{eq:power_law}) that
\begin{align}\label{eq:k_term_error}
\| \theta(f) - \theta_k(f) \|_{\ell^2_M} \le C M^{-1/2} k^{-r}
\end{align}
where $r \deq 1/p - 1/2$, and $C$ is some constant that depends only on $p$. Hence, given any $f \in \cF_d(p)$ and denoting by $f_k$ the function obtained from it by truncating all but the $k$ largest Fourier--Walsh coefficients, we get from Parseval's identity that
\begin{equation}\label{eq:approx_error}
\| f - f_k \|_{L^2(\mu_d)} \le C M^{-1/2} k^{-r}.
\end{equation}
Thus, the assumption that $f$ belongs to the sparsity class $\cF_d(p)$ for some $p$ can be interpreted qualitatively as saying that the behavior of $f$ is strongly influenced by a small number of subsets of the $d$ covariates. The number of these influential subsets decreases as $p \to 0$.

\section{Density estimation via recursive Walsh thresholding}
\label{sec:recursive}

Let $X_1,\ldots,X_n$ be independent random variables in $\bin^d$ with common unknown density $f$. We wish to estimate $f$ on the basis of this sample. For densities defined on the Euclidean space, nonparametric estimators based on hard or soft thresholding of empirically estimated coefficients of the target density in a suitably chosen basis (e.g.,~a wavelet basis) attain near-minimax rates of convergence of the squared-error risk over rich classes of densities \cite{DJKP96,HalKerPic98}. Thresholding is a means of controlling the bias-variance trade-off.

Several authors have investigated the use of term-by-term thresholding rules for density estimation on the binary hypercube. There, one begins by computing the empirical estimates
\begin{equation}
\wh{\theta}_s = \frac{1}{n}\sum^n_{i=1} \chr_s(X_i)
\label{eq:empirical_coeffs}
\end{equation}
of the Fourier--Walsh coefficients of $f$, and then forming the thresholded estimator
\begin{equation}\label{eq:thresh}
	\wh{f} = \sum_{s \in \bin^d} I_{\{T(\wh{\theta}_s) \ge \lambda_n\}} \wh{\theta}_s \chr_s,
\end{equation}
where $T(\cdot)$ is some real-valued statistic and $I_{\{\cdot\}}$ is the indicator function. Based on the observation that
\begin{equation}\label{eq:coeff_var}
 \var \wh{\theta}_s = \frac{1}{n}\left(\frac{1}{M} - \theta^2_s\right),
\end{equation}
while the squared bias incurred by omitting the term
$\wh{\theta}_s\chr_s$ from the estimator \eqref{eq:thresh} is $\theta^2_s$, Ott and Kronmal \cite{OttKro76} considered the ideal thresholded estimator
\begin{equation}\label{eq:ideal_estimator}
\wh{f}^* = \sum_{s \in \bin^d} I_{\{ \theta^2_s > 1/M(n+1)\}}\wh{\theta}_s \chr_s.
\end{equation}
Clearly, $\wh{f}^*$ is impractical because the thresholding criterion depends on the unknown coefficients $\theta_s$. Instead, Ott and Kronmal \cite{OttKro76} suggested that one could mimic the ideal estimator (\ref{eq:ideal_estimator}) by replacing $\theta^2_s$ in the thresholding criterion by the unbiased estimator
$(n\wh{\theta}^2_s - 1/M)/(n-1)$, leading to the practical estimator
\begin{equation}\label{eq:WT_estimator}
\wh{f}_{\rm WT} = \sum_{s \in \bin^d} I_{\{ \wh{\theta}^2_s > 2/M(n+1)\}} \wh{\theta}_s \chr_s,
\end{equation}
where WT stands for ``Walsh thresholding." This estimator was further improved by Liang and Krishnaiah \cite{LiaKri85} and Chen, Krishnaiah and Liang \cite{CheKriLia89}, who replaced the hard thresholding rule in (\ref{eq:WT_estimator}) with shrinkage rules.

The main disadvantage of such termwise thresholding is the need to compute empirical estimates of all $M = 2^d$ Fourier--Walsh coefficients. While this is not an issue when $d$ is comparable to $\log n$, it is clearly impractical when $d \gg \log n$. In order to alleviate this difficulty, we will consider a recursive thresholding approach, which will allow us to reject whole {\em groups} of empirical coefficients based on efficiently implementable thresholding rules. The main idea behind this approach is motivated by the following argument.

Given some $1 \le k \le d$, we can represent any function $f \in L^2(\mu_d)$ with the Fourier--Walsh coefficients $\{\theta_s : s \in \bin^d\}$ as 
\begin{align*}
f &= \sum_{u \in \bin^k} \sum_{v \in \bin^{d-k}} \theta_{uv} \chr_{uv} \\
&= \sum_{u \in \bin^k} \Bigg(\sum_{v \in \bin^{d-k}} \theta_{uv} \chr_v \Bigg) \otimes \chr_u \\
&\equiv \sum_{u \in \bin^k} f_u \otimes \chr_u,
\end{align*}
where, for each $u \in \bin^k$, $f_u \deq \sum_{v \in \bin^{d-k}} \theta_{uv} \chr_v$ is a function in $L^2(\mu_{d-k})$. The Fourier--Walsh coefficients of $f_u$ are precisely those coefficients of $f$ that are indexed by $s \in \bin^d$ with $\pi_k(s) = u$. By Parseval's identity, we have
$$
W_u \deq \| f_u \|^2_{L^2(\mu_{d-k})} = \sum_{v \in \bin^{d-k}} \theta_{uv}^2.
$$
This leads to the following observation: for any $\lambda > 0$,
$$
W_u < \lambda \mbox{ for some } u \in \bin^k \quad \Rightarrow \quad \theta^2_{uv} < \lambda \mbox{ for every } v \in \bin^{d-k}.
$$
The usefulness of this observation for our purposes comes from the
fact that we can represent the strings $s \in \bin^d$, and hence the
elements of the Walsh system in $L^2(\mu_d)$, by the leaves of a
complete binary tree of depth $d$. Suppose we wanted to pick out only
those coefficients of $f$ whose squared magnitude exceeds some
threshold $\lambda$. If we knew that $W_u \le \lambda$ for some $u \in
\bin^k$, then this would tell us that the square of every coefficient
corresponding to a leaf descending from $u$ does not exceed
$\lambda$. Hence, we could start at the root of the tree and at each
internal node $u$ that has not yet been visited check whether $W_u \ge
\lambda$; if not, then we would delete $u$ and all of its descendants
from the tree without having to compute explicitly the corresponding
  coefficients. At the end of the process (i.e., when we get to the
leaves), we will be left only with those $s \in \bin^d$ for which
$\theta^2_s \ge \lambda$. If $f \in \cF_d(p)$ for some $p$, then the
resulting squared $L^2$ error will be
$$
\sum_{s \in \bin^d} I_{\{\theta^2_s < \lambda\}} \theta^2_s \le CM^{-1}(M\lambda)^{-2r/(2r+1)},
$$
where $r = 1/p - 1/2$, as before.

We will follow this reasoning in constructing our density estimator. 
We begin by developing a suitable estimator for $W_u$. To do that, we shall rely on the following lemma (see Appendix~\ref{app:fu} for the proof):

\begin{lemma}\label{lm:fu} For any density $f$ on $\bin^d$, any $1 \le k \le d$, and any $u \in \bin^k$, we have
$$
f_u(y) = \E_f \left\{ \chr_u(\pi_k(X)) I_{\{\sigma_k(X) = y\}}\right\}, \forall y \in \bin^{d-k}
$$
and
$$
W_u = \E_f \left\{ \chr_u(\pi_k(X)) f_u(\sigma_k(X)) \right\}.
$$
\end{lemma}

\noindent This lemma suggests that, for each $1 \le k \le d$ and each $u \in \bin^k$, an empirical estimate of $W_u$ can be obtained by
\begin{eqnarray}
\wh{W}_u &=& \frac{1}{n}\sum^n_{i=1}\chr_u(\pi_k(X_i))\left[\frac{1}{n}\sum^n_{j=1}\chr_u(\pi_k(X_j))I_{\{\sigma_k(X_i) = \sigma_k(X_j)\}} \right] \nonumber\\
&=&\frac{1}{n^2} \sum^n_{i=1} \sum^n_{j=1} \chr_u(\pi_k(X_i))\chr_u(\pi_k(X_j))I_{\{\sigma_k(X_i) = \sigma_k(X_j)\}}.
\label{eq:weight_estimate}
\end{eqnarray}
Although this is a biased estimator, it has the following useful property (see Appendix~\ref{app:wu} for the proof):

\begin{lemma}\label{lm:wu} For any $1 \le k \le d$ and any $u \in \bin^k$,
\begin{equation}
\wh{W}_u = \sum_{v \in \bin^{d-k}} \wh{\theta}^2_{uv},
\label{eq:weight_estimate_direct}
\end{equation}
where each $\wh{\theta}_{uv}$ is an empirical estimate of $\theta_{uv}$ computed according to (\ref{eq:empirical_coeffs}).
\end{lemma}

\noindent Another advantage of computing $\wh{W}_u$ indirectly via (\ref{eq:weight_estimate}), rather than (\ref{eq:weight_estimate_direct}), is that, while the latter requires $O(2^{d-k}n)$ operations, the former requires only $O(n^2d)$ operations. This can amount to significant computational savings when $k < d - \log (nd)$. When $k \ge d - \log (nd)$, it becomes more efficient to use the direct estimator (\ref{eq:weight_estimate_direct}).

Now that we have a way of estimating $W_u$, we can define our density estimation procedure. Provided the threshold scales appropriately with the sample size, we will be able to achieve a good balance between the estimation error (variance) and the approximation error (squared bias) and attain near-minimax rates of convergence. In our analysis, we shall actually consider a more flexible strategy: for every $1 \le k \le d$, we shall compare the estimate $\widehat{W}_u$ of $W_u$ to a threshold that depends not only on the sample size $n$, but also on $k$. More specifically, we will let
\begin{equation}
\lambda_{k,n} = \frac{\alpha_k}{n}, \qquad 1 \le k \le d
\label{eq:thresholds}
\end{equation}
where the sequence $\{\alpha_k\}^d_{k=1}$ satisfies $\alpha_1 \ge \alpha_k \ge \ldots \ge \alpha_d > 0$. In particular, this set-up covers the following two extreme cases:
\begin{enumerate}
\item $\alpha_k = {\rm const}$ for all $k$ -- this covers the standard case of always comparing against the same threshold (that depends on $n$)
\item $\alpha_k = {\rm const} \cdot 2^{d-k}$ -- this corresponds to thresholding not the sum of (a particular subset of) the coefficients, but their {\em average}.
\end{enumerate}
As we shall see, this $k$-dependent scaling will allow us to flexibly trade off the expected $L^2$ error and the computational complexity of the resulting estimator. Now we describe our density estimator. Given the sequence $\boldsymbol{\lambda} = \{\lambda_{k,n}\}^d_{k=1}$, define the set
\begin{equation}\label{eq:accepting_set}
A(\boldsymbol{\lambda}) \deq \left\{ s \in \bin^d : \wh{W}_{\pi_k(s)} \ge \lambda_{k,n}, \forall 1 \le k \le d \right\}
\end{equation}
and consider the density estimate
\begin{equation}\label{eq:RWT}
\wh{f}_{\rm RWT} \deq \sum_{s \in \bin^d} I_{\{ s \in A(\boldsymbol{\lambda}) \}} \wh{\theta}_s \chr_s,
\end{equation}
where RWT stands for ``recursive Walsh thresholding." To implement this estimator on a computer, we call the routine {\sc RecursiveWalsh}, shown as Algorithm~\ref{alg:RWT}, with $u = \varnothing$ (the empty string, corresponding to the root of the tree) and with the desired threshold sequence $\boldsymbol{\lambda}$. The factors of $1/2$ in the updates for $\wh{W}_{u_0}$ and $\wh{W}_{u_1}$ arise because of the factorization property \eqref{eq:walsh_factorize} of the Walsh basis functions: for any $k \ge 0$ and any $s,x,x' \in \bin^{k+1}$ we have
\begin{align*}
	\chr_u(x)\chr_u(x') = \chr_{\pi_k(s)}(\pi_k(x))\chr_{\pi_k(s)}(\pi_k(x'))\left(-\frac{1}{\sqrt{2}}\right)^{s^{(k)}(x^{(k)}+x'^{(k)})}.
\end{align*}

\begin{algorithm}[h]
\caption{{\sc RecursiveWalsh}$(u,\boldsymbol{\lambda})$}

\begin{algorithmic}\label{alg:RWT}
\STATE $k \leftarrow {\rm length}(u)$
\IF{$k = d$}
	\STATE $\wh{\theta}_u \leftarrow \frac{1}{n}\Sum^n_{i=1} \chr_u(X_i)$
	\IF{$\wh{\theta}^2_u \ge \lambda_{d,n}$}
		\STATE {\bf output} $u,\wh{\theta}_u$
	\ENDIF
	\STATE {\bf return}
\ENDIF

\STATE $u_0 \leftarrow 0u$

\STATE $u_1 \leftarrow 1u$

\STATE $\wh{W}_{u_0} \leftarrow \frac{1}{2n^2} \Sum^n_{i=1}\Sum^n_{j=1} \chr_{u}(\pi_{k}(X_i))\chr_{u}(\pi_
{k}(X_j))I_{\{\sigma_{k+1}(X_i) = \sigma_{k+1}(X_j)\}}$

\IF{$\wh{W}_{u_0} \le \lambda_{k+1,n}$}
	\STATE {\bf return}
\ELSE
	\STATE {\sc RecursiveWalsh}$(u_0,\boldsymbol{\lambda})$
\ENDIF

\STATE $\wh{W}_{u_1} \leftarrow \frac{1}{2n^2} \Sum^n_{i=1}\Sum^n_{j=1} \chr_{u}(\pi_{k}(X_i))\chr_{u}(\pi_{k}(X_j))(-1)^{X^{(k+1)}_i + X^{(k+1)}_j} I_{\{\sigma_{k+1}(X_i) = \sigma_{k+1}(X_j)\}}$

\IF{$\wh{W}_{u_1} \le \lambda_{k+1,n}$}
	\STATE {\bf return}
\ELSE
	\STATE {\sc RecursiveWalsh}$(u_1,\boldsymbol{\lambda})$
\ENDIF
\end{algorithmic}
\end{algorithm}

\subsection{Analysis of performance}

Let us denote by $\cF^{+,1}_d(p)$ the set of all probability densities in $\cF_d(p)$:
\begin{equation}\label{eq:den}
	\cF^{+,1}_{d}(p) \deq \left\{ f \in \cF_d(p): f \ge 0, \sum_{x \in \bin^d}f(x) = 1\right\}.
\end{equation}
Our first main result is that, with appropriately tuned thresholds, the estimator (\ref{eq:RWT}) adapts to unknown sparsity of the Fourier--Walsh representation of $f$:

\begin{theorem}\label{thm:RWT} There exist absolute constants $C_1,C_2 > 0$, such that the following holds. Suppose the threshold sequence $\boldsymbol{\lambda} = \{\lambda_k\}^d_{k=1}$ is chosen in such a way that $\alpha_k \ge C_1 d/M$ for all $k$, where $M \equiv 2^d$. Then for all $n \le 2 r M^{2r+1}$ the estimator (\ref{eq:RWT}) satisfies
\begin{equation}\label{eq:l2_risk}
\sup_{f \in \cF^{+,1}_d(p)} \E_f \| f - \wh{f}_{\rm RWT} \|^2_{L^2(\mu_d)} \le \frac{C_2}{M} \left(\frac{\log M}{n}\right)^{2r/(2r+1)}
\end{equation}
for all $0 < p \le 1$, where, as before, $r = 1/p - 1/2$. Moreover, if $\alpha_k \ge C_1 d \log n / M$ for all $k$, then the risk of \eqref{eq:RWT} is bounded as
\begin{equation}\label{eq:l2_risk_adaptive}
	\sup_{f \in \cF^{+,1}_d(p)} \E_f \| f - \wh{f}_{\rm RWT} \|^2_{L^2(\mu_d)} \le \frac{C_2}{M} \left(\frac{\log M \log n}{n}\right)^{2r/(2r+1)}
\end{equation}
for all $n$.
\end{theorem}

\begin{remark} Positivity and normalization issues. {\em As is the case with orthogonal series estimators, $\wh{f}_{\rm RWT}$ may not necessarily be a bona fide density. In particular, there may be some $x \in \bin^d$ such that $\wh{f}_{\rm RWT}(x) < 0$, and it may happen that $\int \wh{f}_{\rm RWT}d\mu_d \neq 1$. In principle, this can be handled by clipping the negative values at zero and renormalizing; this procedure can only improve the expected $L^2$ error. In practice this may be computationally expensive when $d$ is very large. If the estimate is suitably sparse, however, the renormalization can be carried out approximately using Monte-Carlo estimates of the appropriate sums. Moreover, in many applications the scaling of the density is not important.\hfill $\square$}
\end{remark}

\begin{remark} Logarithmic factors in the risk bound. {\em For each $0 < p \le 1$, the bound \eqref{eq:l2_risk} will hold for all $n \le 2r M^{2r+1}$. Since
$$
2r M^{2r+1} = \left(\frac{2}{p}-1\right)M^{2/p} \ge M^{2/p}, \qquad 0 < p \le 1
$$
it follows from Theorem~\ref{thm:minimax} below that $\wh{f}_{\rm RWT}$ with thresholds $\lambda_1 = \ldots = \lambda_d \sim nd/M$ is minimax-optimal for each $0 < p \le 1$, assuming that the sample sample size $n$ satisfies $d \equiv \log M \preceq n \preceq M^{2/p}$. For very small sample sizes $n \preceq \log M$ and for very large sample sizes $n \ge M^{2/p}$, $\wh{f}_{\rm RWT}$ will be suboptimal. With a more conservative choice of thresholds, $\lambda_1 = \ldots = \lambda_d \sim nd\log n/M$, the bound \eqref{eq:l2_risk_adaptive} will hold for all values of $n$. In particular, in this case the minimax rate will be attained, up to logarithmic factors, for all values of $p \in (0,1)$ simultaneously in the moderate sample regime $\log M \preceq n \preceq M^{2/p}$.  \hfill $\square$}
\end{remark}

\noindent Our second main result is a lower bound on the minimax $L^2$ risk attainable by any estimator over $\cF^{+,1}_d(p)$. It shows, in particular, that our recursive estimator $\wh{f}_{\rm RWT}$ is minimax for ``moderate'' sample sizes $\log M \preceq n \preceq M^{2/p}$. For large sample sizes, $n \succeq M^{2/p}$, $\wh{f}_{\rm RWT}$ is no longer optimal --- in particular, it is outperformed by both the simple histogram estimator
$$
\wh{f}_{\rm hist}(x) = \frac{1}{n}\sum^n_{i=1}I_{\{X_i = x\}}
$$
and by the unthresholded orthogonal series estimator
$$
\wh{f}(x) = \sum_{s \in \bin^d} \wh{\theta}_s \chr_s(x),
$$
both of which attain the optimal $O(1/n)$ risk. The precise statement is as follows:

\begin{theorem}\label{thm:minimax} Consider the problem of estimating
  an unknown $f \in \cF^{+,1}_d(p)$ from $n$ i.i.d.\ samples
  $X_1,\ldots,X_n$. Then the following statements hold:
\begin{enumerate}
	\item Suppose that $\log M \le n \le M^{2(1-\epsilon)/p}$ for some $\epsilon \in (0,1)$. Then there exists a positive constant $C = C(p,\epsilon)$, such that
	\begin{align}\label{eq:small_sample_LB}
		\inf_{\wh{f}_n}\sup_{f \in \cF^{+,1}_d(p)} \E_f \| \wh{f}_n - f \|^2_{L^2(\mu_d)} \ge \frac{C}{M}\left(\frac{\log M}{n}\right)^{\frac{2r}{2r+1}}.
	\end{align}
where, as before, $M = 2^d$.
\item Suppose that $n \ge M^{2/p}$ and $M\ge 4$. Then there exists an absolute constant $C > 0$, such that
\begin{align}\label{eq:large_sample_LB}
	\inf_{\wh{f}_n}\sup_{f \in \cF^{+,1}_d(p)} \E_f \| \wh{f}_n - f \|^2_{L^2(\mu_d)} \ge \frac{C}{n}.
\end{align}
\end{enumerate}
\end{theorem}

\noindent Our third, and final, main result bounds the running time of the algorithm used for computing $\wh{f}_{\rm RWT}$:

\begin{theorem}\label{thm:comp} Fix any $f \in \cF_d(p)$. Given any $\delta \in (0,1)$, provided each $\alpha_k$ is chosen so that
\begin{equation}
	nC_1 (2^k a^2_{k,n} \wedge 2^{k/2}a_{k,n}) \ge \frac{\log (2^k d/\delta)}{\log e},
	\label{eq:alpha_choice}
\end{equation}	
where
$$
a_{k,n} \deq \frac{1}{5}\sqrt{\frac{\alpha_k}{n}} - \sqrt{\frac{C^2_2}{2^k n}}
$$
and $C_1,C_2 > 0$ are certain absolute constants, then Algorithm~\ref{alg:RWT} runs in
\begin{equation}
O\left(n^2d \left(\frac{n}{M}\right)^{p/2} K(\boldsymbol{\alpha},p)\right)
\label{eq:time_complexity}
\end{equation}
time with probability at least $1-\delta$, where $K(\boldsymbol{\alpha},p) \deq \sum^d_{k=1} \alpha^{-p/2}_k$ and, as before, $M = 2^d$.
\end{theorem}

\begin{remark} Trade-off between time complexity and MSE. {\em By controlling the rate at which the sequence $\alpha_k$ decays with $k$, we can trade off MSE against complexity. Consider the following two extreme cases: (1) $\alpha_1 = \ldots = \alpha_d \sim 1/M$ and (2) $\alpha_k \sim 2^{d-k}/M$. The first case, which reduces to the term-by-term thresholding, achieves the same bias-variance trade-off as the Ott--Kronmal estimator \cite{OttKro76}. However, it has $K(\boldsymbol{\alpha},p) = O(M^{p/2}d)$, resulting in $O(d^2n^{2+p/2})$ complexity. The second case, which leads to a very severe estimator that will tend to reject a lot of coefficients, has MSE of $O(n^{-2r/(2r+1)}M^{-1/(2r+1)})$, but $K(\boldsymbol{\alpha},p) = O(M^{p/2})$, leading to a considerably better $O(dn^{2+p/2})$ complexity. From the computational viewpoint, it is preferable to use rapidly decaying thresholds. However, this reduction in complexity will be offset by a corresponding increase in MSE. In fact, using RWT with exponentially decaying $\alpha_k$'s in practice is not advisable as its low complexity is mainly due to the fact that it will tend to reject even the big coefficients very early on, especially when $d$ is large. To achieve a good balance between complexity and MSE, a moderately decaying threshold sequence might be best, e.g.,~$\alpha_k \sim (d-k+1)^m/M$ for some $m \ge 
1$. As $p \to 0$, the effect of $\boldsymbol{\lambda}$ on complexity becomes negligible, and the complexity tends to $O(n^2d)$. \hfill $\square$}
\end{remark}

\begin{remark} Incoherence. {\em Note that for any of the above choices of
    $\alpha_k$, the proposed method requires ${\rm polylog}(M)$
    operations. One intuitive explanation for why such fast
    computation is possible is that the Walsh basis is
    ``incoherent'' (to use term common in compressed sensing and group
    testing literature) with the canonical basis of $L^2(\mu_d)$. Similar
    polylog computation results were achieved by Gilbert {\em et al.} in the context of fast sparse Fourier approximation \cite{gilbertMansour,GilbertSPIE} and group testing \cite{GilbertStrauss}. Their strategies also had
    connections to the Goldreich--Levin algorithm \cite{GolLev89}, as well as to the work of Kushilevitz and Mansour on sparse Boolean function estimation \cite{KusMan93,Man94}. \hfill $\square$}
\end{remark}

\section{Proofs of the theorems}
\label{sec:proofs}

In this section we prove our three main results. However, before proceeding to the proofs, we collect all the technical tools that we will be using: moment bounds, concentration inequalities, and an approximation-theoretic lemma pertaining to class $\cF^{+,1}_d(p)$.

\subsection{Preliminaries}
\label{ssec:prelims}

\subsubsection{Moment bound} We will need the following result of Rosenthal \cite{Ros72}. Let $U_1,\ldots,U_n$ be i.i.d.~random variables with $\E U_i = 0$ and $\E U^2_i \le \sigma^2$. Then for any $m \ge 2$ there exists some $c_m$ such that
\begin{equation}\label{eq:rosenthal}
\E \left|n^{-1}\sum^n_{i=1} U_i\right|^m \le c_m \left(\frac{\sigma^m}{n^{m/2}} + \frac{\E |U_1|^m}{n^{m-1}}\right).
\end{equation}

\subsubsection{Concentration bounds} We will need the well-known Hoeffding inequality: if $U_1,\ldots,U_n$ are i.i.d.\ bounded random variables such that $\E U_i = 0$ and $|U_i| \le b < \infty$ for all $1 \le i \le n$, then
\begin{equation}\label{eq:hoeffding}
	\Pr \left( \left|\frac{1}{n}\sum^n_{i=1}U_i \right| > t \right) \le 2 \exp\left(-\frac{nt^2}{2b^2}\right)
\end{equation}
The following result is due to Talagrand \cite{Tal94}. Let $U_1,\ldots,U_n$ be i.i.d.~random variables, let $\varepsilon_1,\ldots,\varepsilon_n$ be independent Rademacher random variables [i.e.,~$\Pr(\varepsilon_i = -1) = \Pr(\varepsilon_i = 1) = 1/2$] also independent of $U_1,\ldots,U_n$, and let $\cF$ be a class of functions uniformly bounded by $L > 0$. Then if there exist some $v,H > 0$ such that $\sup_{g \in \cF} \var g(U) \le v$ and
\begin{equation}\label{eq:rademacher}
\E \left\{ \sup_{g \in \cF} \sum^n_{i=1} \varepsilon_i g(U_i) \right\} \le nH
\end{equation}
for all $n$, then there are universal constants $C_1$ and $C_2$ such that, for every $\tau > 0$,
\begin{equation}\label{eq:talagrand}
\Pr \left(\sup_{g \in \cF} \nu_n(g) \ge \tau + C_2H \right) \le \exp \left\{ - nC_1 \left(\frac{\tau^2}{v} \wedge \frac{\tau}{L}\right)\right\},
\end{equation}
where
\begin{equation}\label{eq:emp}
\nu_n(g) \deq \frac{1}{n} \sum^n_{i=1} g(U_i) - \E g(U), \qquad \forall g \in \cF
\end{equation}
is the empirical process indexed by $\cF$.

\begin{remark}{\em Typically, some additional regularity conditions on $\cF$ are needed to ensure measurability of the supremum $\sup_{g \in \cF}\nu_n(g)$ of the empirical process \eqref{eq:emp}. However, when $U$ takes values in a finite set, as is the case in this paper, there is no need for such conditions because any uniformly bounded class of real-valued functions on a finite set is separable: it contains a countable subset $\cF_0$, such that for any $g \in \cF$ there exists a sequence $g_1,g_2,\ldots \in \cF_0$ converging to $g$ pointwise. Such a separability property ensures measurability of suprema over $\cF$ \cite[p.~110]{Vaart_Wellner}.}\end{remark}

\subsubsection{Large separated subsets of $\cF^{+,1}_d(p)$}

In the sequel, we will be interested in large subsets of the class of densities $\cF^{+,1}_d(p) \subset \cF_d(p)$, whose elements are separated from one another by a given fixed amount, as measured by the norm $\| \cdot \|_{L^2(\mu_d)}$. The following lemma, whose proof is given in Appendix~\ref{app:proofs}, will be useful:

\begin{lemma}\label{lm:separated} Let $r = 1/p - 1/2$. Let $s_1,\ldots,s_M$, $M=2^d$, be the lexicographic ordering of the elements of $\bin^d$. Given a positive real parameter $a$ and an integer $k \in \{1,\ldots,M-1\}$, let $\Theta(M,k,a) \subset \R^{M-1}$ consist of $(M-1)$-dimensional real vectors having exactly $k$ nonzero components, each of which is equal to either $a$ or $-a$. With this, define the set $\cF(k,a) \subset L^2(\mu_d)$ by
\begin{align*}
	\cF(k,a) = \left\{ f : \theta_{s_1}(f) = \frac{1}{\sqrt{M}}, \,\, \left(\theta_{s_j}(f)\right)_{2 \le j \le M} \in \Theta(M,k,a) \right\}.
\end{align*}
Suppose that $k$ and $a$ are such that
\begin{align}\label{eq:ka_conditions}
	ka \le \frac{1}{2\sqrt{M}} \qquad \text{and} \qquad a \le \frac{1}{\sqrt{M}}(k+1)^{-1/p}.
\end{align}
Then the following statements hold:
\begin{enumerate}
	\item The set $\cF(k,a)$ is contained in $\cF^{+,1}_d(p)$. 
	\item For any two $f,f' \in \cF(k,a)$ we have
	\begin{align}\label{eq:KL_bound}
		D(f \| f') \le 2 M k a^2,
	\end{align}
	where $D(\cdot \| \cdot)$ is the Kullback--Leibler divergence (relative entropy) \cite{CoverThomas}.
	\item If $k=M-1$, then there exists a set $\tilde{\cF}(M-1,a) \subset \cF(M-1,a)$ with the following properties:
	\begin{itemize}
		\item $\| f - f' \|^2_{L^2(\mu_d)} \ge (M-1)a^2$ for all $f,f' \in \tilde{\cF}(M-1,a)$ with $f \neq f'$
		\item $\log |\tilde{\cF}(M-1,a)| \ge (M-1)/8$
	\end{itemize}
	\item If $M-1 \ge 4k$, then there exists a set $\tilde{\cF}(k,a) \subset \cF(k,a)$ with the following properties:
	\begin{itemize}
		\item $\| f - f' \|^2_{L^2(\mu_d)} \ge ka^2$ for all $f,f' \in \tilde{\cF}(k,a)$ with $f \neq f'$
		\item $\log |\tilde{\cF}(k,a)| \ge 0.233\, k \left(\log \frac{M-1}{k} + 1\right)$
	\end{itemize}
\end{enumerate}
\end{lemma}

\subsection{Proof of Theorem~\ref{thm:RWT}}

Let us decompose the squared $L^2$ error as
\begin{align*}
 \| f - \wh{f}_{\rm RWT} \|^2_{L^2(\mu_d)} &= \sum_s I_{\{s \in A(\boldsymbol{\lambda})\}} (\theta_s - \wh{\theta}_s)^2 + \sum_s I_{\{s \in A(\boldsymbol{\lambda})^c\}} \theta^2_s \equiv   T_1 + T_2.
\end{align*}
We start by observing that any $s \in A(\boldsymbol{\lambda})$ necessarily satisfies $\wh{W}_s \equiv \wh{\theta}^2_s \ge \lambda_{d,n}$, while for any $s \in A(\boldsymbol{\lambda})^c$ there exists some $1 \le k \le d$ such that $\wh{W}_{\pi_k(s)} < \lambda_{k,n}$, which implies, in turn, that $\wh{\theta}^2_{\pi_k(s)t} < \lambda_{k,n}$ for all $t \in \bin^{d-k}$ and, in particular, that $\wh{\theta}^2_s < \lambda_{k,n} \le \lambda_{1,n}$. Therefore, defining the sets
$$
A_1 = \{ s \in \bin^d : \wh{\theta}^2_s \ge \lambda_{d,n} \} \quad \mbox{and} \quad A_2 = \{ s \in \bin^d : \wh{\theta}^2_s < \lambda_{1,n} \},
$$
we can bound $T_1$ and $T_2$ as
$$
T_1 \le \sum_s I_{\{ s \in A_1 \}} (\theta_s - \wh{\theta}_s)^2
$$
and
$$
T_2 \le \sum_s I_{\{ s \in A_2 \}} \theta^2_s.
$$
Further, defining
$$
B = \Big\{ s \in \bin^d : \theta^2_s < \lambda_{d,n}/2 \Big\} \quad \mbox{and} \quad  S = \Big\{ s \in \bin^d : \theta^2_s \ge 3\lambda_{1,n}/2 \Big\},
$$
we can write
\begin{align*}
T_1 &\le \sum_s I_{\{ s \in A_1 \cap B\}} (\theta_s - \wh{\theta}_s)^2 + \sum_s I_{\{ s \in A_1 \cap B^c \}} (\theta_s - \wh{\theta}_s)^2 \equiv T_{11} + T_{12},
\end{align*}
and
\begin{align*}
T_2 &\le \sum_s I_{\{s \in A_2 \cap S\}} \theta^2_s + \sum_s I_{\{s \in A_2  \cap S^c\}} \theta^2_s \equiv T_{21} + T_{22}.
\end{align*}
Applying (\ref{eq:concentration}) and \eqref{eq:coeff_var}, we get
\begin{align}
 \E T_{12} & \le \sum_{s \in B^c}  \E (\theta_s - \wh{\theta}_s)^2  \nonumber \\
 & \le \frac{1}{Mn}\left|\left\{ s  : \theta^2_s \ge \lambda_{d,n}/2\right\}\right| \nonumber  \\
 &\le \frac{1}{Mn}\left( \frac{2}{M\lambda_{d,n}}\right)^{p/2} \nonumber \\
 & = \frac{1}{Mn} \left( \frac{2 n}{M \alpha_d}\right)^{p/2} \nonumber \\
 & = \frac{1}{M} n^{p/2 - 1} \underbrace{\left(\frac{2}{2^d \alpha_d }\right)^{p/2}}_{\le 1} \nonumber \\
 & \le \frac{1}{M}n^{-2r/(2r+1)}. \label{eq:T12}
 \end{align}
To bound $T_{22}$ we apply (\ref{eq:approx_error}):
 \begin{align}
 \E T_{22} \le \sum_{s \in \bin^d} I_{\{ \theta_s^2 < 3\alpha_1/2n\}} \theta^2_s &\le \frac{C}{M}\left(\frac{M \alpha_1}{n}\right)^{2r/(2r+1)}\nonumber\\
&\le \frac{C}{M}\left(\frac{\log M}{n}\right)^{2r/(2r+1)}.\label{eq:T22}
\end{align}
In order to deal with the large-deviation terms $T_{11}$ and $T_{21}$, we will need some moment and concentration bounds which are listed in Section~\ref{ssec:prelims}. First, using Cauchy--Schwarz, we get
\begin{equation}\label{eq:T11_1}
\E T_{11} \le \sum_s \sqrt{\E (\theta_s - \wh{\theta}_s)^4  \cdot \Pr( A_1 \cap B) }
\end{equation}
To estimate the fourth moment in (\ref{eq:T11_1}), we apply the bound (\ref{eq:rosenthal}) to $U_i = \chr_s(X_i) - \theta_s$, $1 \le i \le n$, and $m = 4$. Then
$$
\E U^2_i = \E (\chr_s(X_i) - \theta_s)^2 = \frac{1}{M} - \theta^2_s \le \frac{1}{M}
$$
and
$$
\E |U_1|^4 \le \frac{1}{M^2},
$$
so that
$$
\E (\theta_s - \wh{\theta}_s)^4 \le c_4 \left(\frac{1}{M^2 n^2} + \frac{1}{M^2 n^3}\right) \le \frac{2c_4}{M^2 n^2}.
$$
To handle the probability of $A_1 \cap B$, we first estimate
$$
|\wh{\theta}_s - \theta_s|^2 = (\theta_s - \wh{\theta}_s)^2 =
\theta^2_s - 2\theta_s\wh{\theta}_s + \wh{\theta}^2_s \ge \theta^2_s -
2|\theta_s\wh{\theta}_s| + \wh{\theta}^2_s = (|\theta_s| -
|\wh{\theta}_s|)^2.
$$
From this we conclude that $\wh{\theta}^2_s \ge \lambda_{d,n}$ and $\theta^2_s < \lambda_{d,n}/2$ together imply
$$
| \wh{\theta}_s - \theta_s | \ge \frac{1}{5}\sqrt{\lambda_{d,n}} = \frac{1}{5} \sqrt{\frac{\alpha_d}{n}}
$$
(the factor of 1/5 is simply a lower bound on $1-1/\sqrt{2}$). Therefore,
$$
\Pr(A_1 \cap B) \le \Pr \left(| \wh{\theta}_s - \theta_s | \ge \frac{1}{5} \sqrt{\frac{\alpha_d}{n}}\right).
$$
Applying Hoeffding's inequality (\ref{eq:hoeffding}) to $U_i = \chr_s(X_i) - \theta_s$, $1 \le i \le n$, with $b = 2/\sqrt{M}$ and using the fact that $\alpha_d \ge C_1 d/M$, we get
\begin{align*}
\Pr\left(| \wh{\theta}_s - \theta_s | \ge \frac{1}{5} \sqrt{\frac{\alpha_d}{n}}\right) &\le 2\exp\left(-{Cd}\right) \le \frac{2}{M^C}.
\end{align*}
for some absolute constant $C > 0$. If $C_1$ is chosen so that $C \ge 2$, then we will have
\begin{equation}
\E T_{11} \le \frac{2\sqrt{ c_4}}{Mn}.\label{eq:T11}
\end{equation}
Finally,
$$
\E T_{21} \le \sum_s \Pr(A_2 \cap S) \theta^2_s.
$$
Using the same argument as above, we can write
$$
\Pr ( A_2 \cap S ) \le \frac{2}{M^2}
$$
(with the same choice of the constant $C_1$ as before). Then
\begin{align*}
	\E T_{21} &\le \frac{2}{M^2}\sum_s \theta^2_s \\
	&\le \frac{2}{M^2} \cdot \frac{1}{M}\sum^M_{m=1}m^{-2/p} \\
	&\le \frac{2}{M^2} \cdot \frac{M^{-(2r+1)}}{2r},
\end{align*}
where the first inequality uses the fact that $f \in \cF_d(p)$, while the second inequality follows from
\begin{align*}
	\sum^M_{m=1} m^{-2/p} &\le \int^M_0 t^{-2/p} dt = \frac{M^{-(2/p-1)}}{2/p-1} = \frac{M^{-2r}}{2r}.
\end{align*}
Since $n \le 2r M^{2r+1}$, $M^{-(2r+1)}/2r \le n^{-1}$. Consequently, we will have
\begin{equation}\label{eq:T21}
\E T_{21} \le \frac{C}{M^2n}.
\end{equation}
Putting together Eqs.~(\ref{eq:T12}), (\ref{eq:T22}), (\ref{eq:T11}), and (\ref{eq:T21}), we get (\ref{eq:l2_risk}). The second bound \eqref{eq:l2_risk_adaptive} is proved along the same lines, except that the extra $\log n$ factor in the thresholds will give $\E T_{11} \le C/Mn$ and $\E T_{21} \le C/Mn$.

\subsection{Proof of Theorem~\ref{thm:minimax}}

The proof of the first part uses a popular information-theoretic technique due to Yang and Barron \cite{YangBarronTechRep,YangBarron}; we only outline the main steps. The first step is to lower-bound the minimax risk by the minimum probability of error in a multiple hypothesis test. Let $\cF_0$ be an arbitrary subset of $\cF^{+,1}_d(p)$. Then
$$
\inf_{\wh{f}}\sup_{f \in \cF^{+,1}_d(p)} \E_f \| \wh{f} - f \|^2_{L^2(\mu_d)} \ge \inf_{\wh{f}}\sup_{f \in \cF_0} \E_f \| \wh{f} - f \|^2_{L^2(\mu_d)}.
$$
In particular, suppose that the set $\cF_0$ is finite, $\cF_0 = \{f^{(1)},\ldots,f^{(N)}\}$, and $\delta$-separated in $L^2(\mu_d)$, i.e.,
\begin{align}\label{eq:separated}
\| f^{(i)} - f^{(j)} \|_{L^2(\mu_d)} \ge \delta, \qquad \forall i,j \in \{1,\ldots,N\}; i \neq j
\end{align}
Then a standard argument \cite{YangBarron} gives
\begin{equation}\label{eq:YB_step1}
	\inf_{\wh{f}} \sup_{f \in \cF_0} \E_f \| \wh{f} - f \|^2_{L^2(\mu_d)} \ge \frac{\delta^2}{4} \min_{\tilde{f}}\Pr \left( \tilde{f} \neq f^{(Z)}\right),
\end{equation}
where the random variable $Z$ is uniformly distributed over the set $\{1,\ldots,N\}$, and the minimum is over all estimators $\tilde{f}$ based on $X^n$ that take values in the packing set $\{f^{(1)},\ldots,f^{(N)}\}$. Applying Fano's inequality \cite{CoverThomas}, we can write
\begin{equation}\label{eq:YB_step2}
\min_{\tilde{f}}\Pr \left( \tilde{f} \neq f^{(Z)}\right) \ge 1 - \frac{I(Z; X^n) + \log 2}{\log N},
\end{equation}
where $I(Z; X^n)$ is the Shannon mutual information \cite{CoverThomas} between the random index $Z \in \{1,\ldots,N\}$ and the observations $X_1,\ldots,X_n \stackrel{\text{i.i.d.}}{\sim} f^{(Z)}$. In this particular case, we have
\begin{align}\label{eq:MI}
	I(Z; X^n) &= \frac{n}{N}\sum^N_{k=1}\sum_{x \in \bin^d} f^{(k)}(x) \log \frac{f^{(k)}(x)}{N^{-1}\sum^N_{\ell=1} f^{(\ell)}(x)} \nonumber\\
	&= \frac{n}{N}\sum^N_{k=1}D(f^{(k)}\| \bar{f}),
\end{align}
where $\bar{f}$ denotes the mixture density $N^{-1}\sum^N_{\ell=1}f^{(\ell)}$. The next step consists in upper-bounding this mutual information. To that end, suppose that there exists some $\Delta > 0$, such that
\begin{align}\label{eq:KL_upper_bound}
	D(f \| f') \le \Delta, \qquad \forall f,f' \in \cF_0.
\end{align}
Using convexity of the relative entropy and \eqref{eq:KL_upper_bound}, for every $k \in \{1,\ldots,N\}$ we have
\begin{align*}
	D(f^{(k)} \| \bar{f}) \le \frac{1}{N}\sum^N_{\ell=1}D(f^{(k)} \| f^{(\ell)}) \le \Delta.
\end{align*}
Substituting this into \eqref{eq:MI}, we see that $I(Z; X^n) \le n\Delta$. Combining this bound with \eqref{eq:YB_step1}, we get
\begin{align}\label{eq:YB_bound}
			\inf_{\wh{f}} \sup_{f \in \cF_0} \E_f \| \wh{f} - f \|^2_{L^2(\mu_d)} \ge \frac{\delta^2}{4}\left[ 1 - \frac{n\Delta + \log 2}{\log N}\right] 
\end{align}
In particular, let $k \in \{1,\ldots,M-1\}$ and $a > 0$ satisfy the conditions \eqref{eq:ka_conditions} of Lemma~\ref{lm:separated}, as well as
\begin{align}\label{eq:ka_conditions_2}
	a^2 \le \frac{C}{Mnk} \log |\tilde{\cF}(k,a)|
\end{align}
for a suitable constant $C > 0$, where $\tilde{\cF}(k,a)$ are the subsets of $\cF^{+,1}_d(p)$ described in Lemma~\ref{lm:separated}. If we let $\cF_0 = \tilde{\cF}(k,a)$, then, by Lemma~\ref{lm:separated}, \eqref{eq:KL_upper_bound} holds with $\Delta = \frac{1}{2n}\log \frac{|\tilde{\cF}(k,a)|}{4}$. This, in conjunction with \eqref{eq:YB_bound}, gives
\begin{align}
		\inf_{\wh{f}} \sup_{f \in \cF_0} \E_f \| \wh{f} - f \|^2_{L^2(\mu_d)} \ge C\delta^2(k,a),
\end{align}
where
$$
\delta(k,a) \deq \min\left\{ \| f - f' \|_{L^2(\mu_d)} : f,f' \in \tilde{\cF}(k,a); f \neq f' \right\}
$$
is the minimal $L^2(\mu_d)$-separation between any two distinct elements of $\tilde{\cF}(k,a)$. We can now consider the following cases:
\begin{enumerate}
	\item Suppose that $M^{2/p} \le n$. Then we take $k = M-1$ and $a^2 = \frac{C}{Mn}$. Because in this case
	\begin{align*}
		k \vee (k+1)^{1/p} = M^{1/p} \le \sqrt{n} \le \frac{C}{a\sqrt{M}},
	\end{align*}
and $\log |\tilde{\cF}(M-1,a)| \ge (M-1)/8$, the conditions \eqref{eq:ka_conditions} and \eqref{eq:ka_conditions_2} will be satisfied for a suitable choice of $C$. Moreover, by Lemma~\ref{lm:separated}, we have
	\begin{align*}
		\delta^2(k,a) &= \delta^2(M-1,a) \ge (M-1)a^2 \ge \frac{C}{n}.
	\end{align*}
	Substituting this into \eqref{eq:YB_bound}, we obtain \eqref{eq:large_sample_LB}.
	\item Suppose that $n \le M^{2(1-\epsilon)/p}$ for some $\epsilon \in (0,1)$. Let $k = C\left(\frac{n}{\log M}\right)^{p/2}$ and $a = C'\left(\frac{\log M}{Mn}\right)^{1/2}$. Then
	\begin{align*}
		k \vee (k+1)^{1/p} &\le (2k)^{1/p} \\
		&= (2C)^{1/p} \left(\frac{n}{\log M}\right)^{1/2} \\
		&= \frac{(2C)^{1/p}C'}{a\sqrt{M}}.
	\end{align*}
	If we choose $C$ and $C'$ in such a way that $(2C)^{1/p}C' \le 1/2$, then \eqref{eq:ka_conditions} will be satisfied. Moreover, we must have
	\begin{align}\label{eq:k_vs_M}
		M-1 \ge 4k = 4C \left(\frac{n}{\log M}\right)^{p/2}.
	\end{align}
	With our assumptions on $n$ and $M$, this will hold for all sufficiently large $M$. Next, we check that \eqref{eq:ka_conditions_2} is satisfied. Assuming that \eqref{eq:k_vs_M} holds, Lemma~\ref{lm:separated} implies that
	\begin{align*}
		\frac{1}{Mnk} \log |\tilde{\cF}(k,a)| &\ge \frac{C}{Mn} \log \frac{M-1}{k} \\
		&\ge \frac{C}{Mn} \log \frac{(M-1)^{2/p}\log M}{n}.
	\end{align*}
	Again, using our assumption on $n$ and $M$, as well as the fact that $a^2 =  \frac{C \log M}{Mn}$, we can guarantee that \eqref{eq:ka_conditions_2} holds, with an appropriate choice of $C = C(p,\epsilon)$. By Lemma~\ref{lm:separated}, we will have
	\begin{align*}
		\delta^2(k,a) &\ge ka^2 \ge \frac{C}{M} \left(\frac{\log M}{n}\right)^{1-\frac{p}{2}}  = \frac{C}{M} \left(\frac{\log M}{n}\right)^{\frac{2r}{2r+1}},
	\end{align*}
	and we obtain \eqref{eq:small_sample_LB}.
\end{enumerate}

\subsection{Proof of Theorem~\ref{thm:comp}}

The time complexity of the algorithm is determined by the number of recursive calls made to {\sc RecursiveWalsh}. Recall that, for each $1 \le k \le d$, a recursive call to {\sc RecursiveWalsh} is made for every $u \in \bin^k$ for which $\wh{W}_u \ge \lambda_{k,n}$. Let us say that a recursive call to {\sc RecursiveWalsh}$(u,\boldsymbol{\lambda})$ is {\em correct} if $W_u \ge \lambda_{k,n}/2$. We will show that, with high probability, only the correct recursive calls are made at every $1 \le k \le d$. The probability of making at least one incorrect recursive call is given by
\begin{align*}
&\Pr\left( \bigcup^d_{k=1} \bigcup_{u \in \bin^{k}} \{ \wh{W}_u \ge \lambda_{k,n}, W_u < \lambda_{k,n}/2 \} \right) \\
& \qquad \le \sum^d_{k=1} \sum_{u \in \bin^{k}} 
\Pr\left( \wh{W}_u \ge \lambda_{k,n}, W_u < \lambda_{k,n}/2 \right).
\end{align*}
For a given $u \in \bin^{k}$, let
\begin{align*}
\wh{f}_u \deq \sum_{v \in \bin^{d-k}} \wh{\theta}_{uv} \chr_v.
\end{align*}
Then $\wh{W}_u \ge \lambda_{k,n}$ and $W_u < \lambda_{k,n}/2$ together imply that
\begin{align*}
\| f_u - \wh{f}_u \|^2_{L^2(\mu_{d-k})} &= \sum_{v \in \bin^{d-k}} (\wh{\theta}_{uv} - \theta_{uv})^2 \\
&= \wh{W}_u - 2 \sum_{v \in \bin^{d-k}} \wh{\theta}_{uv} \theta_{uv} + W_u \\
&\ge \wh{W}_u - 2 \sqrt{\wh{W}_u W_u} + W_u \\
&= \left(\sqrt{\wh{W}_u} - \sqrt{W_u}\right)^2 \\
& \ge \left(1-\frac{1}{\sqrt{2}}\right)^2 \lambda_{k,n} \\
&\ge \lambda_{k,n}/25.
\end{align*}
Now, as shown in Appendix~\ref{app:epr}, for each $u \in \bin^{k}$, the norm $\| f_u - \wh{f}_u \|_{L^2(\mu_{d-k})}$  can be expressed as a supremum of an empirical process over a suitable function class, to which we can then apply Talagrand's bound (\ref{eq:talagrand}) with $L = 1/\sqrt{2^k}$, $v = 1/2^k$, and $H = 1/\sqrt{2^kn}$. Hence,
\begin{align*}
\Pr ( \wh{W}_u \ge \lambda_{k,n}, W_u < \lambda_{k,n}/2) &\le \Pr ( \| f_u - \wh{f}_u \|^2_{L^2(\mu_{d-k})} \ge \sqrt{\lambda_{k,n}}/5) \\
&\le \exp\left\{ - nC_1(2^k a^2_{k,n} \wedge 2^{k/2} a_{k,n} )\right\},
\end{align*}
where for each $1 \le k \le d$,
$$
a_{k,n} = \frac{1}{5}\sqrt{\frac{\alpha_k}{n}} - \sqrt{\frac{C^2_2}{2^kn}}.
$$
Here, $C_1$ and $C_2$ are the absolute constants in Talagrand's bound (\ref{eq:talagrand}). Given $\delta > 0$, if we choose $\alpha_k$ according to (\ref{eq:alpha_choice}), then
$$
\Pr ( \wh{W}_u \ge \lambda_{k,n}, W_u < \lambda_{k,n}/2) \le \frac{\delta}{d2^{k}}, \qquad \forall u \in \bin^k.
$$
Summing over $1 \le k \le d$ and $u \in \bin^k$, we see that, with probability at least $1-\delta$, only the correct recursive calls will be made.

Next, we give an upper bound on the number of the correct recursive calls. For each $1 \le k \le d$, $W_u \ge \lambda_{k,n}/2$ implies that there exists at least one $v \in \bin^{d-k}$ such that $\theta^2_{uv} \ge \lambda_{k,n}/2$. Since for every $1 \le k \le d$ each $\theta_s$ contributes to exactly one $W_u$, we have by the pigeonhole principle that
\begin{eqnarray*}
\left| \left\{ u \in \bin^k : W_u \ge \lambda_{k,n}/2 \right\} \right| &\le& \left| \left\{ s \in \bin^d: \theta^2_s \ge \lambda_{k,n}/2 \right\} \right| \\
&\le& \left(\frac{2}{M\lambda_{k,n}}\right)^{p/2},
\end{eqnarray*}
where in the second line we used (\ref{eq:concentration}). Hence, the number of the correct recursive calls in Algorithm~\ref{alg:RWT} is bounded by
$$
N = \sum^d_{k=1}\left(\frac{2}{M\lambda_{k,n}}\right)^{p/2} = \left(\frac{2n}{M}\right)^{p/2}\sum^d_{k=1} \alpha^{-p/2}_k.
$$
At each recursive call, we compute an estimate of the corresponding $W_{u0}$ and $W_{u1}$, which requires $O(n^2d)$ operations. Therefore, with probability at least $1-\delta$, the time complexity of the algorithm is given by (\ref{eq:time_complexity}).

\section{Simulations}
\label{sec:simulation}

Although an extensive empirical evaluation is outside the scope of this paper, we have implemented the proposed estimator, and now present some simulation results to demonstrate its small-sample performance on synthetic data in low- and high-dimensional regimes. In the low-dimensional regime, it is feasible to obtain the ``ground truth" by exhaustively computing  all the $2^d$ Walsh coefficients and to compare it with our estimate. In the high-dimensional regime, our comparison is based on the density values at randomly generated samples. Additionally, we present a number of computational strategies that greatly enhance computational efficiency in the high-dimensional regime.

\subsection{Low-dimensional simulations}

\begin{figure}[h]
\centerline{
\begin{tabular}{cc}
\includegraphics[width=0.5\textwidth]{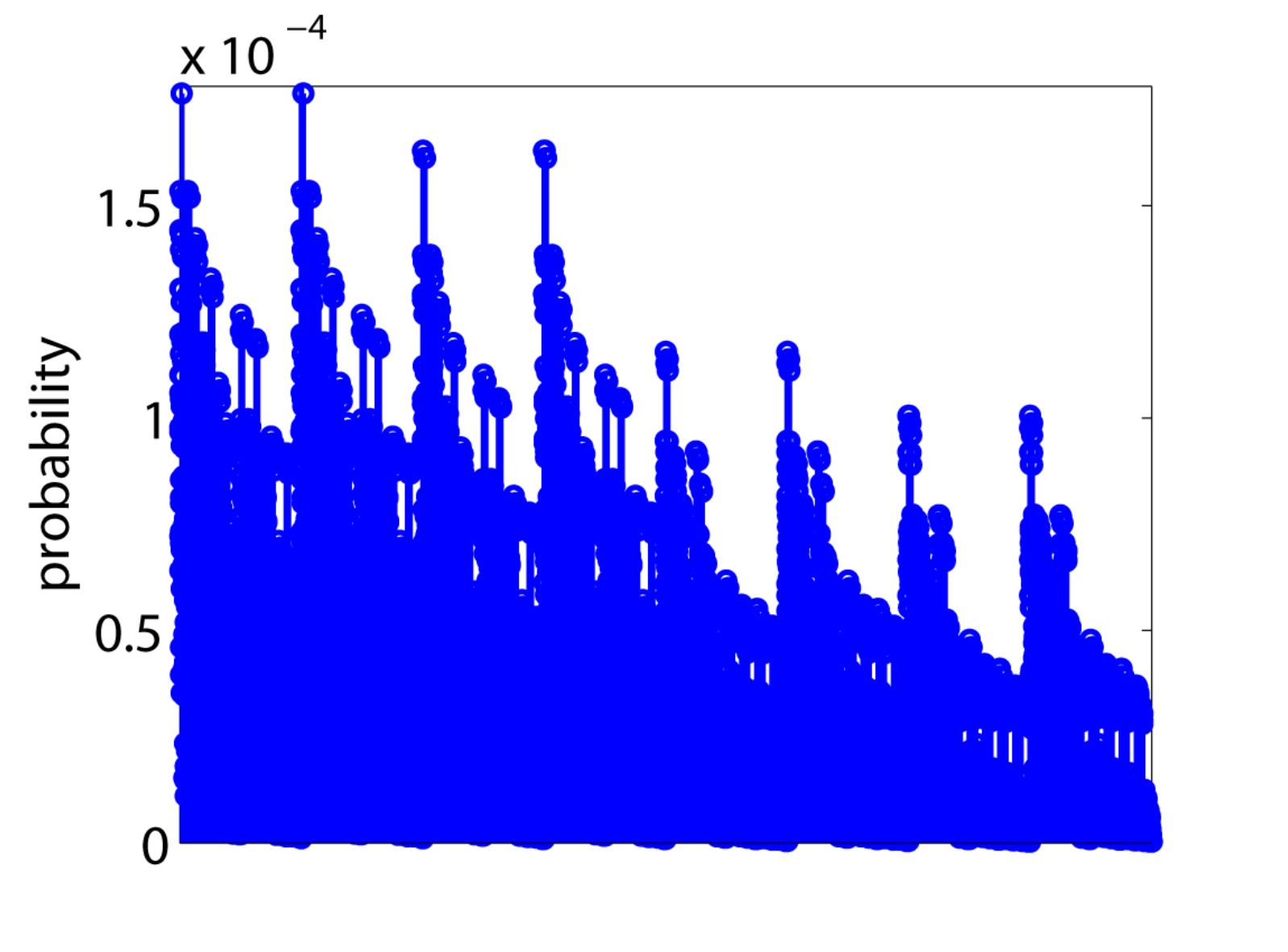} & \includegraphics[width=0.5\textwidth]{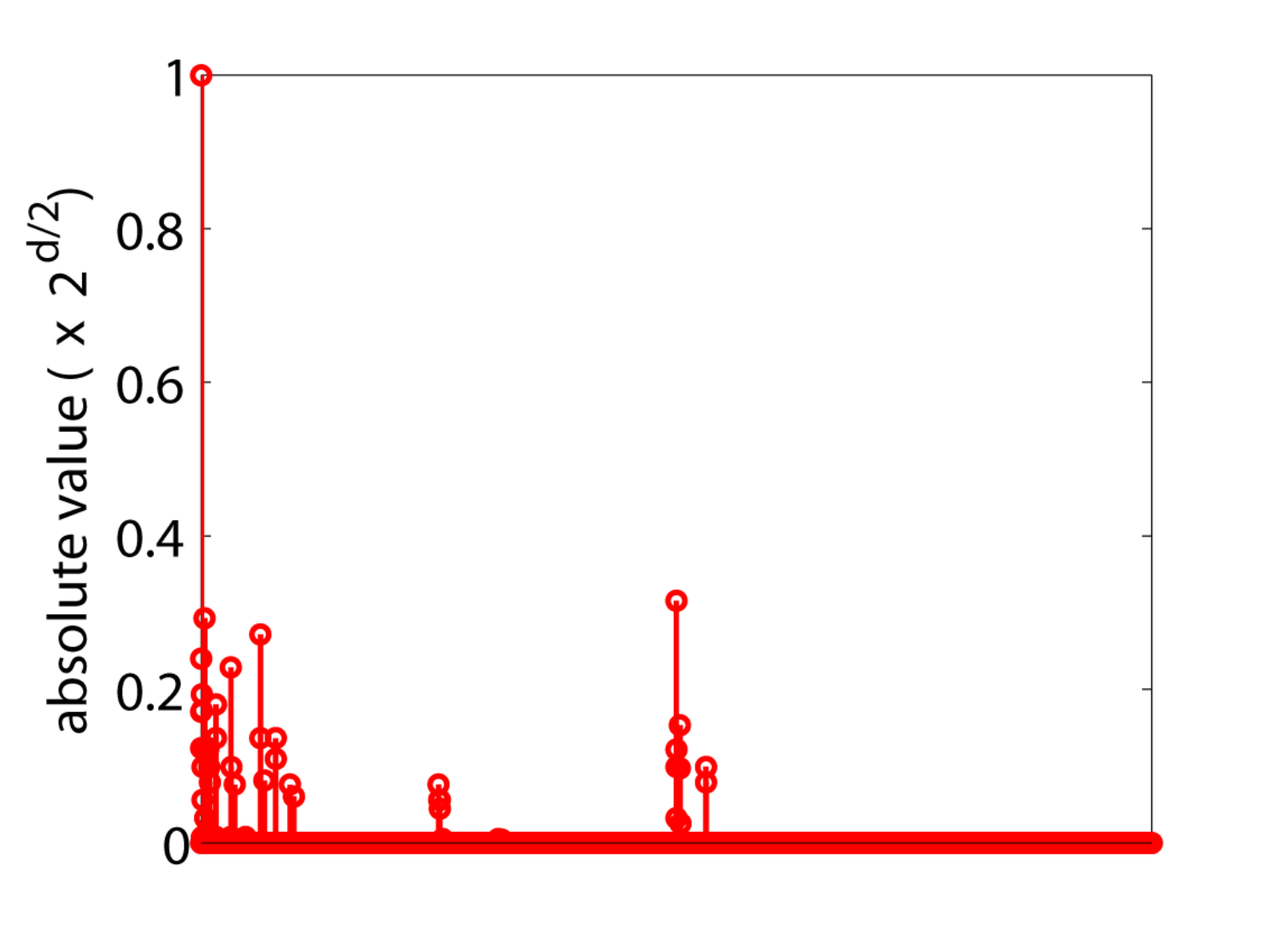} \\
\footnotesize{Ground truth ($f$)} & \\
& \\
\includegraphics[width=0.5\textwidth]{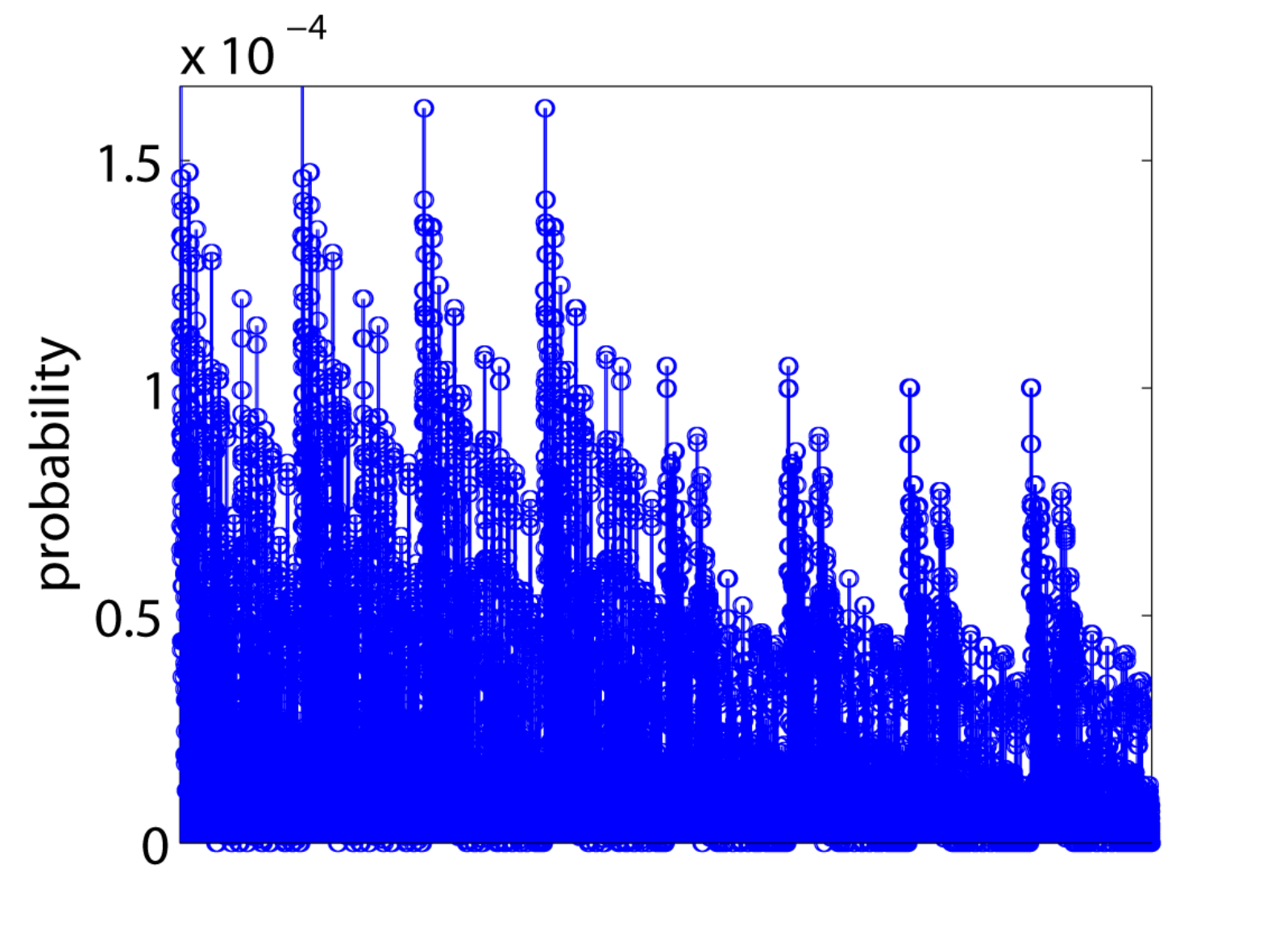} &
\includegraphics[width=0.5\textwidth]{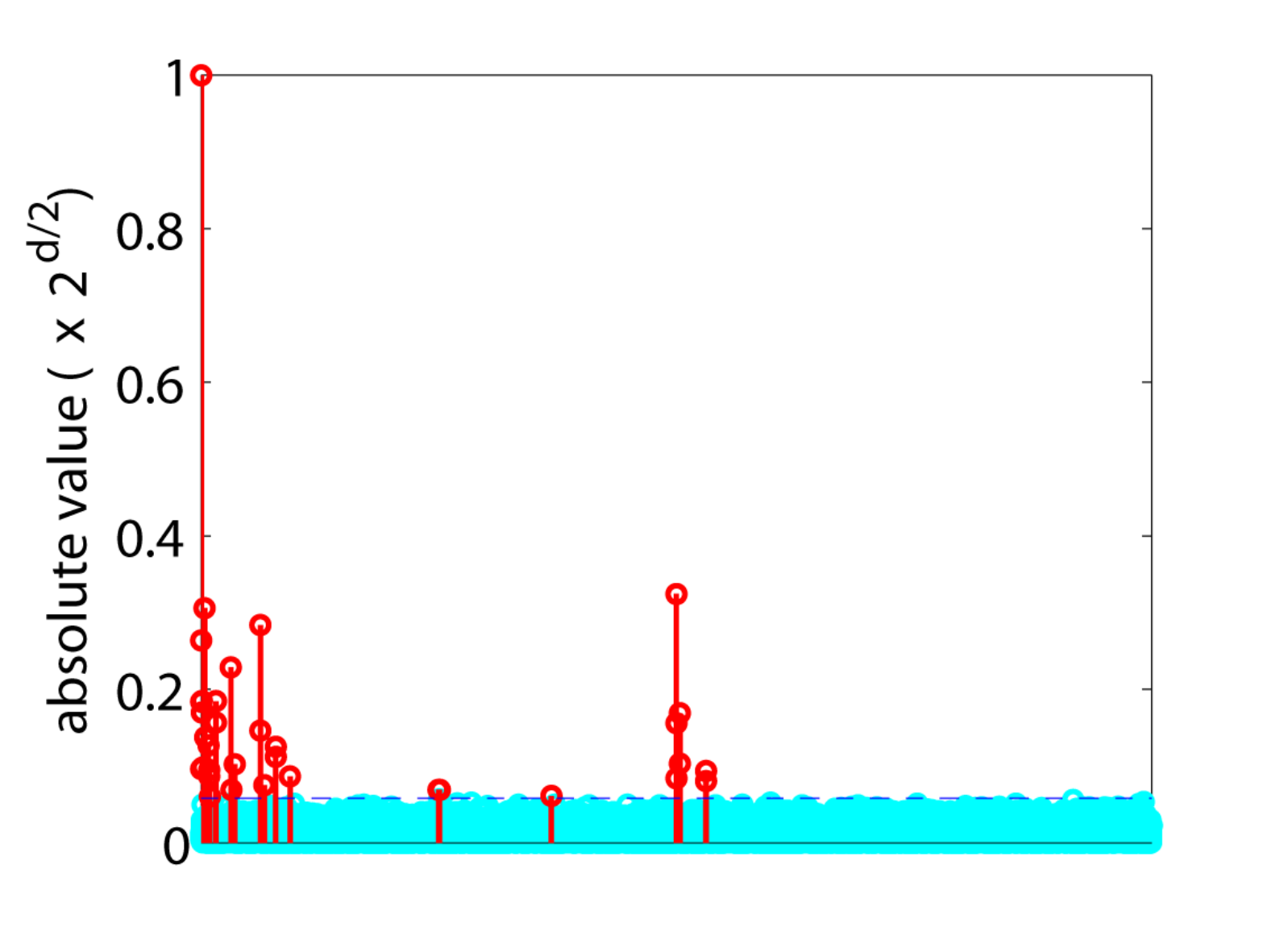} \\
\footnotesize{$\wh{f}_{\rm RWT}$, $n=500$} & \\
& \\
\includegraphics[width=0.5\textwidth]{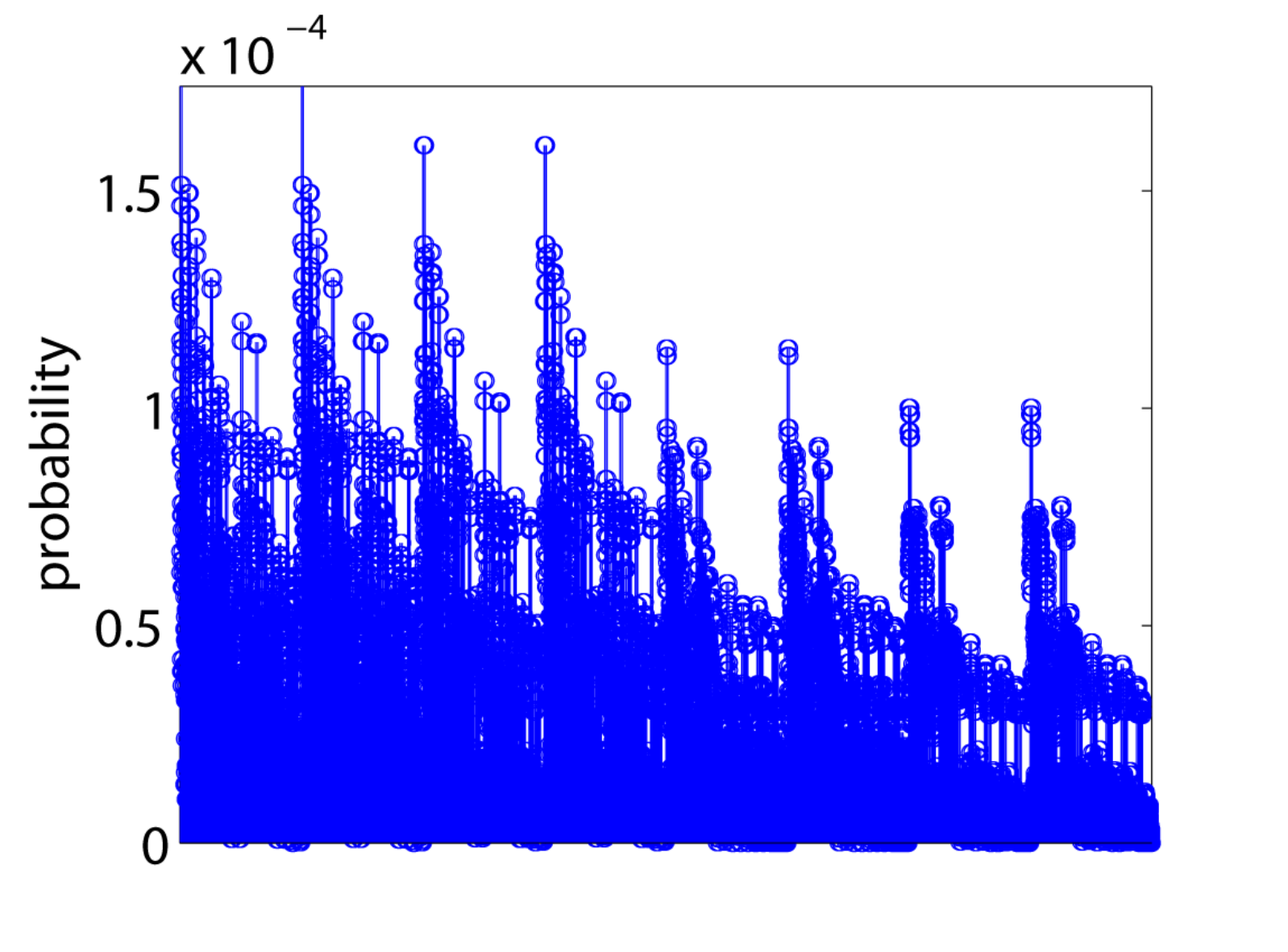} &
\includegraphics[width=0.5\textwidth]{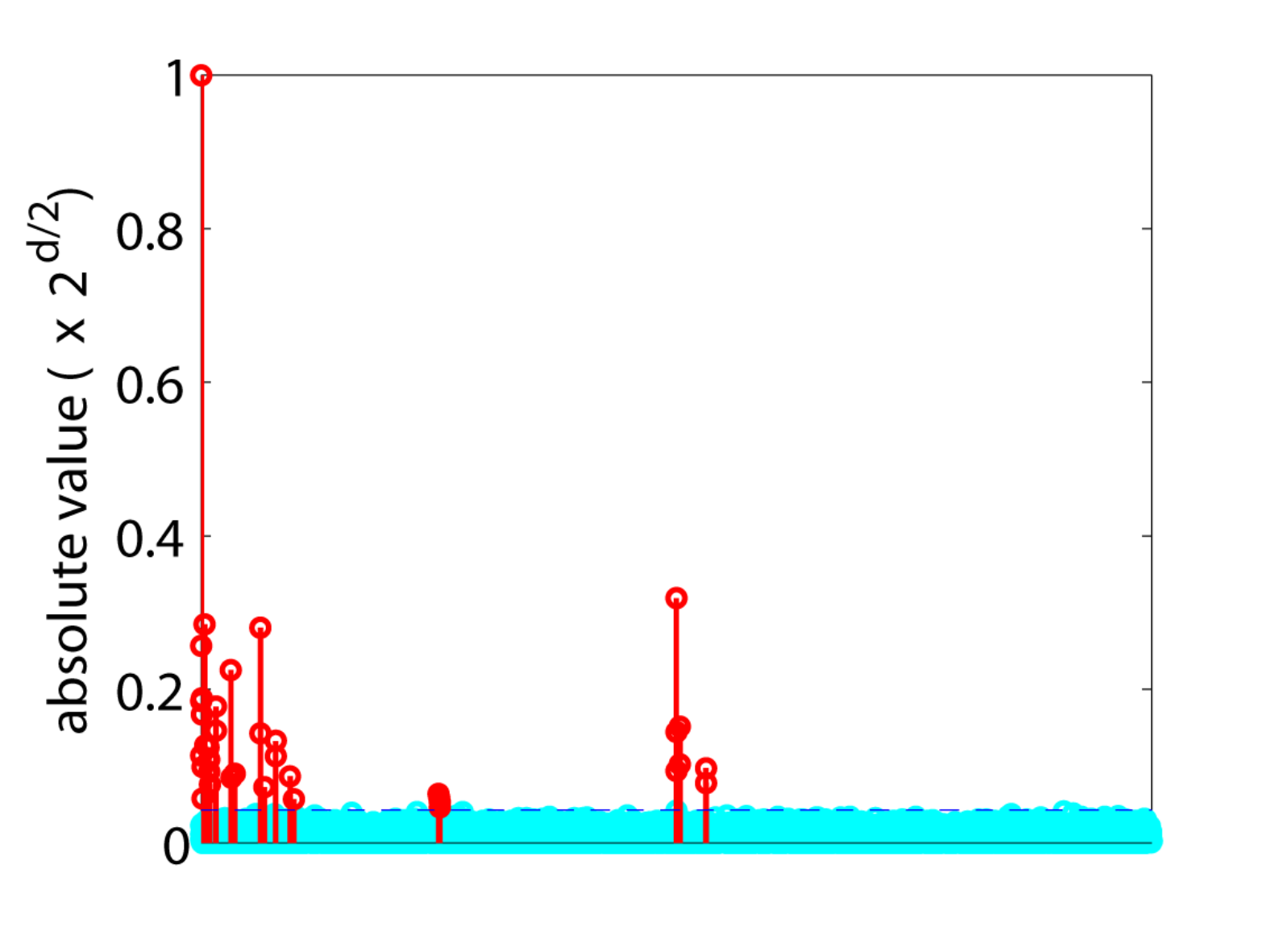} \\
\footnotesize{$\wh{f}_{\rm RWT}$, $n=10000$} &
\end{tabular}}
\caption{\small Ground truth (top) and estimated density for $n=5000$ (middle) and $n=10000$ (bottom) with constant thresholding. Left column: true and estimated probabilities (clipped at zero and renormalized) arranged in lexicographic order. Right column: absolute values of true and estimated Walsh coefficients arranged in lexicographic order. For the estimated densities, the coefficient plots also show the threshold level (dotted line) and absolute values of the rejected coefficients (lighter color). \label{fig:density}}
\end{figure}

We generated synthetic observations from a mixture density $f$ on a $15$-dimensional binary hypercube. The mixture has $10$ components, where each component is a product density with $12$ randomly chosen covariates having ${\rm Bernoulli}(1/2)$ distributions, and the other three having ${\rm Bernoulli}(0.9)$ distributions. For $d=15$, it is still feasible to quickly compute the ground truth, consisting of $32768$ values of $f$ and its Walsh coefficients. These values are shown in Fig.~\ref{fig:density} (left). As can be seen from the coefficient profile in the bottom of the figure, this density is clearly sparse. Fig.~\ref{fig:density} also shows the estimated probabilities and the Walsh coefficients for sample sizes $n=5000$ (middle) and $n=10000$ (right).

\begin{figure}[h]
\centerline{
\setlength{\tabcolsep}{-0.04cm}
\begin{tabular}{cc} 
\includegraphics[width=0.5\textwidth]{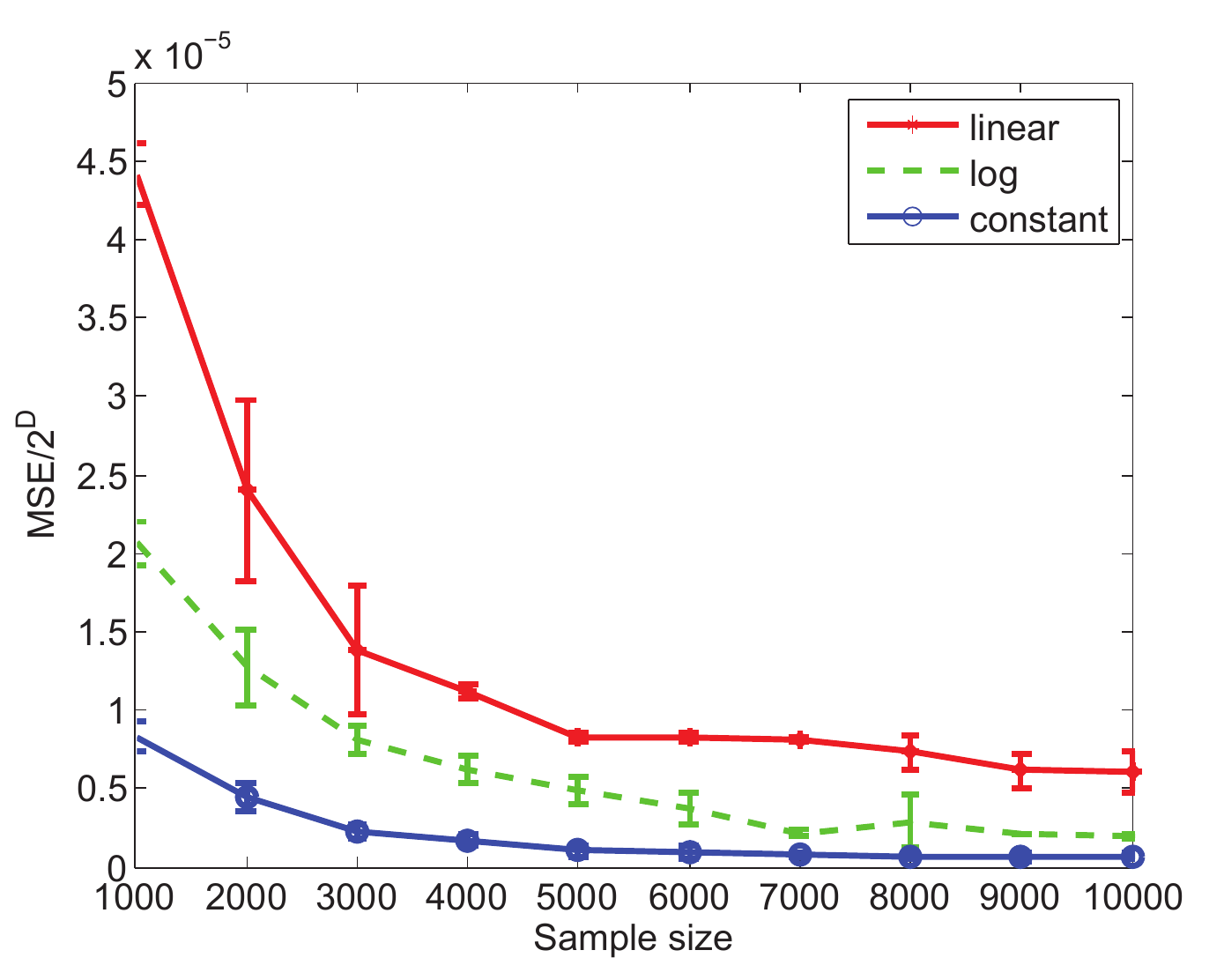} & 
\includegraphics[width=0.5\textwidth]{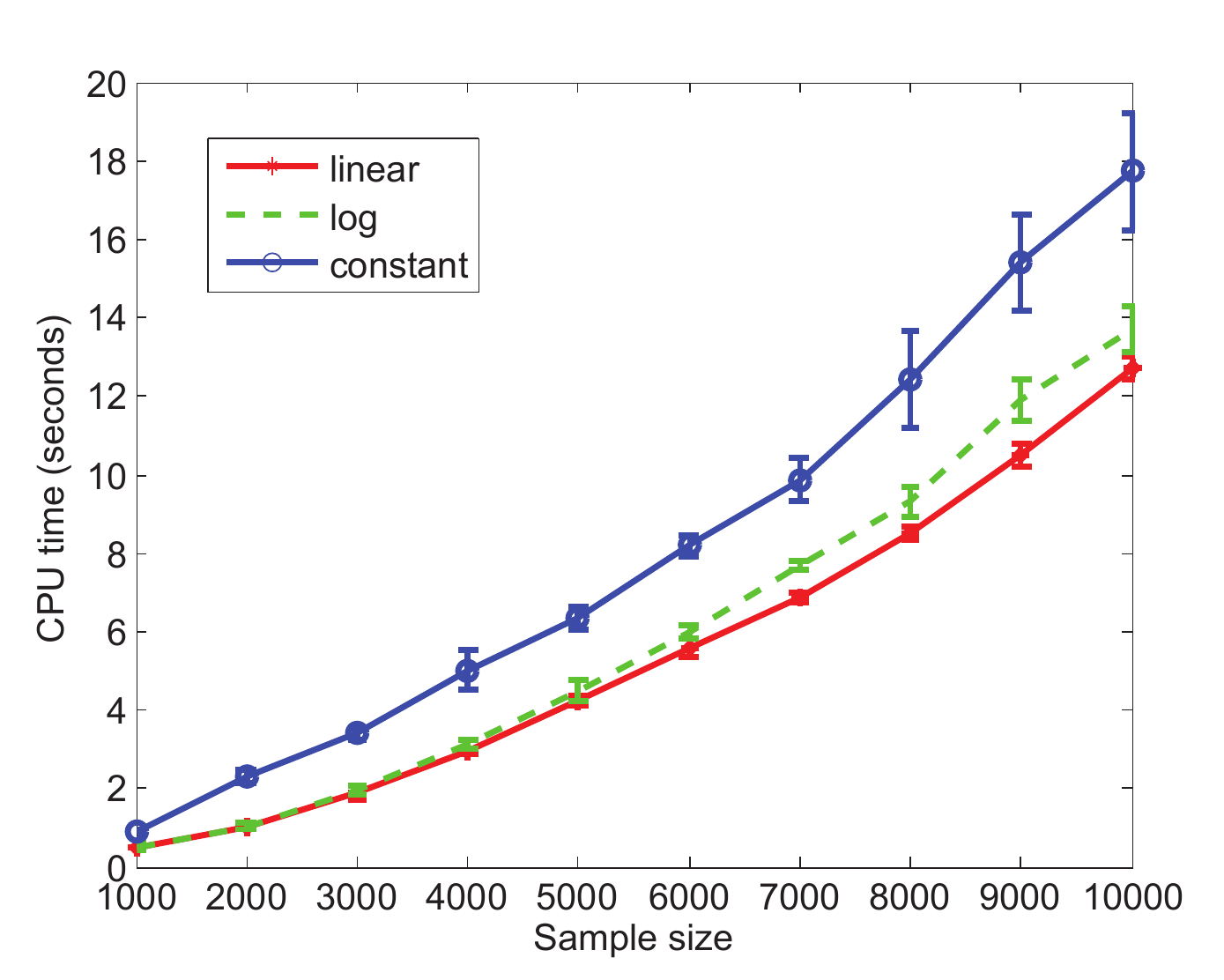} \\
\footnotesize{(a)} & \footnotesize{(b)} \\
\includegraphics[width=0.5\textwidth]{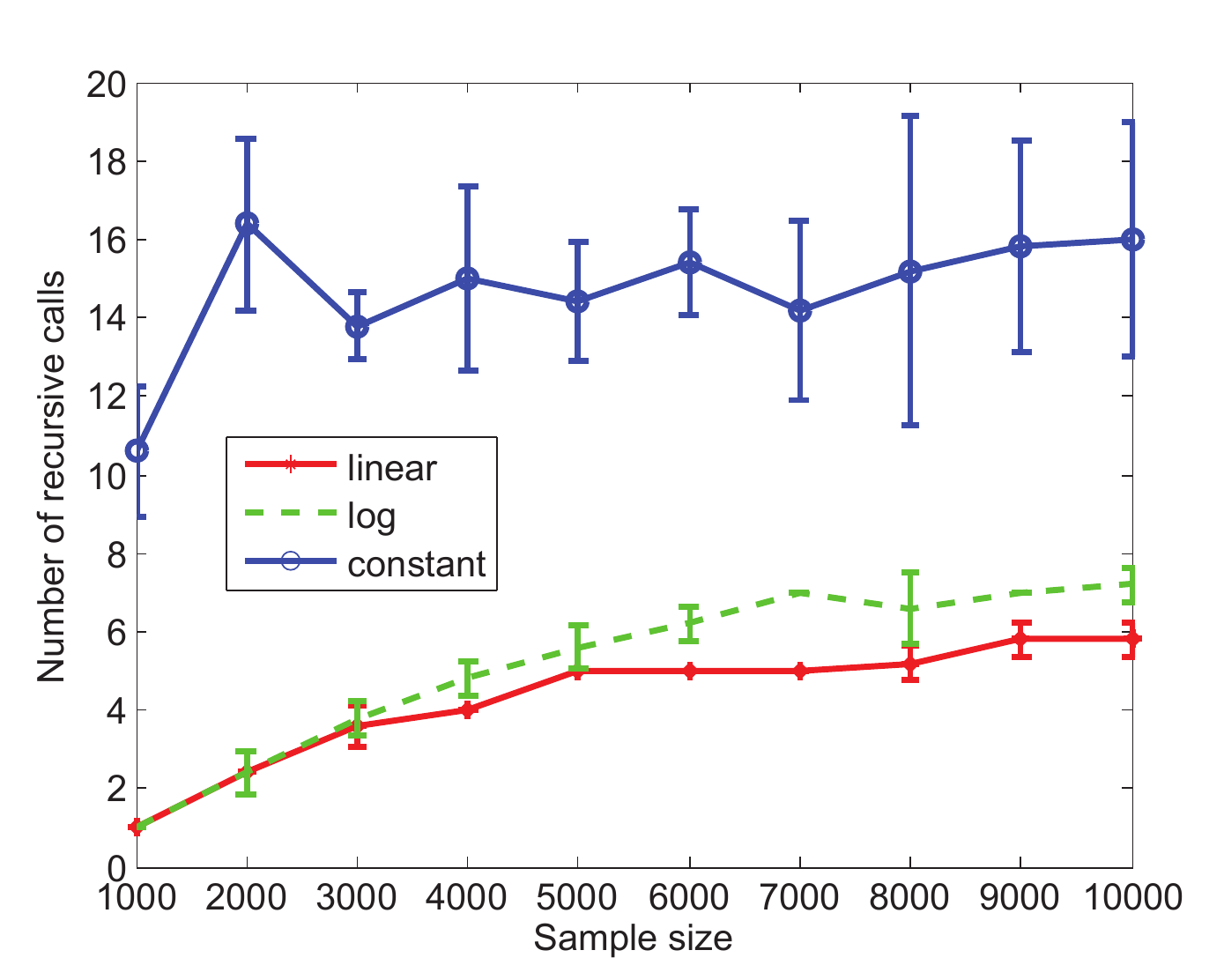} &
\includegraphics[width=0.5\textwidth]{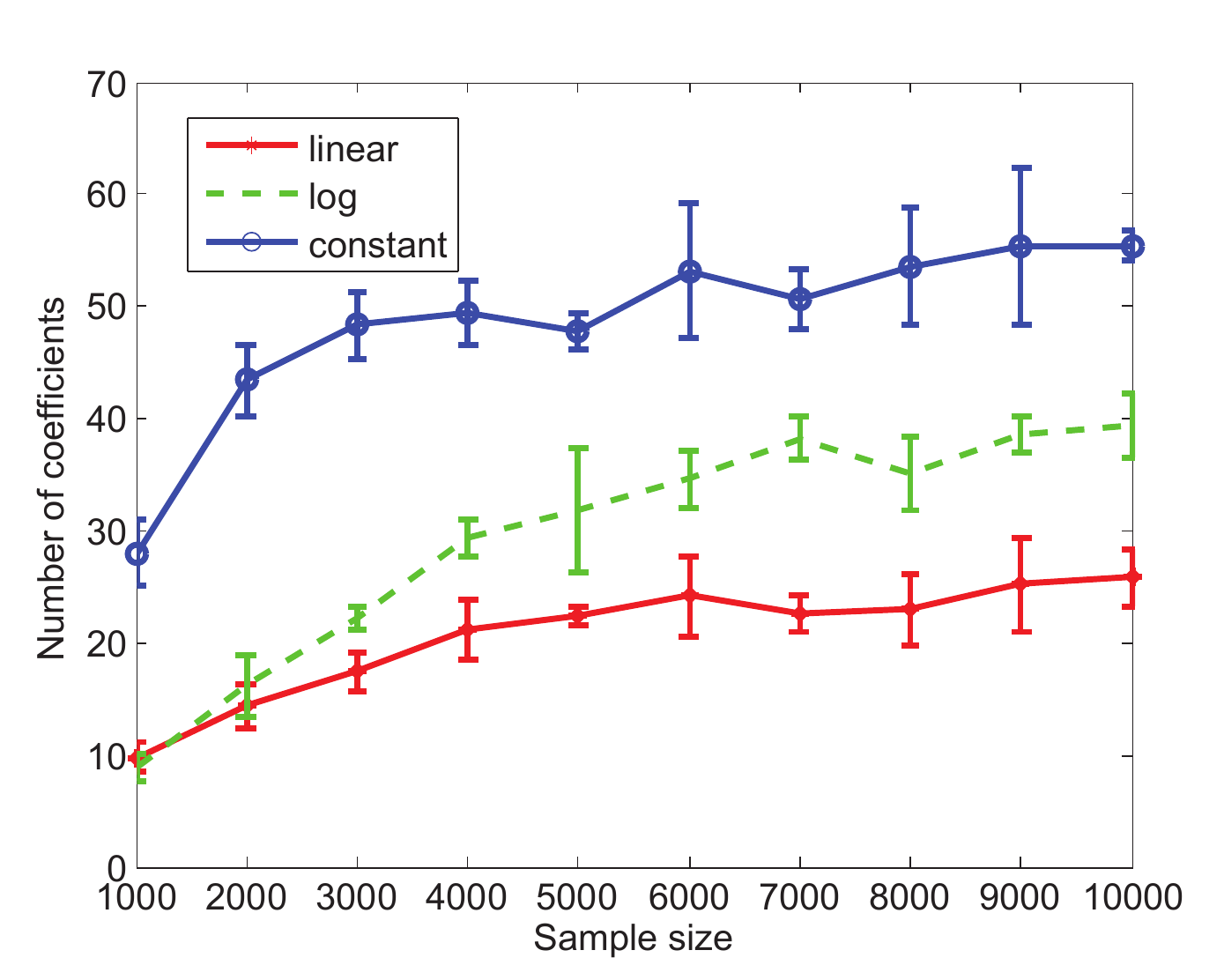} \\
 \footnotesize{(c)} & \footnotesize{(d)}
\end{tabular}}
\caption{\small Small-sample performance of $\wh{f}_{\rm RWT}$ in estimating $f$ with three different thresholding schemes: (a)~MSE; (b)~running time (in seconds); (c)~number of recursive calls; (d)~number of coefficients retained by the algorithm. All results are averaged over five independent runs for each sample size (the error bars show the standard deviations). \label{fig:performance}}
\end{figure}

To study the trade-off between MSE and complexity, we implemented three different thresholding schemes: (1) constant, $\lambda_{k,n} = 2/(2^dn)$, (2) logarithmic, $\lambda_{k,n} = 2 \log (d-k+2)/(2^dn)$, and (3) linear, $\lambda_{k,n} = 2(d-k+1)/(2^dn)$. The thresholds at $k=d$ are set to twice the variance of the empirical estimate of any coefficient whose value is zero; this forces the estimator to reject empirical coefficients whose values cannot be reliably distinguished from zero. Occasionally, spurious coefficients get retained, as can be seen in Fig.~\ref{fig:density} (middle) for the estimate for $n=5000$. Fig~\ref{fig:performance} shows the performance of $\wh{f}_{\rm RWT}$. Fig.~\ref{fig:performance}(a) is a plot of MSE vs.~sample size. In agreement with the theory, MSE is the smallest for the constant thresholding scheme [which is simply an efficient recursive implementation of a term-by-term thresholding estimator with $\lambda_n \sim 1/Mn$], and then it increases for the logarithmic and for the linear schemes. Fig.~\ref{fig:performance}(b,c) shows the running time (in seconds) and the number of recursive calls made to {\sc RecursiveWalsh} vs.~sample size. The number of recursive calls is a platform-independent way of gauging the computational complexity of the algorithm, although it should be kept in mind that each recursive call has $O(n^2d)$ overhead. Also, the number of recursive calls depends on whether a binary or $N$-ary tree is utilized. The $N$-ary tree scheme is explained in detail below, in Section \ref{sec:optim}.  We have used $N=256$ in the simulations, as this setting leads to much reduced computations times vs. a binary tree.

The running time increases polynomially with $n$, and is the largest for the constant scheme, followed by the logarithmic and the linear schemes. We see that, while the MSE of the logarithmic scheme is fairly close to that of the constant scheme, its complexity is considerably lower, in terms of both the number of recursive calls and the running time. In all three cases, the number of recursive calls decreases with $n$ due to the fact that weight estimates become increasingly accurate with $n$, which causes the expected number of false discoveries (i.e.,~making a recursive call at an internal node of the tree only to reject its descendants later) to decrease. Finally, Fig.~\ref{fig:performance}(d) shows the number of coefficients retained in the estimate. This number grows with $n$ as a consequence of the fact that the threshold decreases with $n$, while the number of accurately estimated coefficients increases.

\begin{figure}[h]
\centerline{
\setlength{\tabcolsep}{-0.04cm}
\begin{tabular}{cc} 
\includegraphics[width=0.5\textwidth]{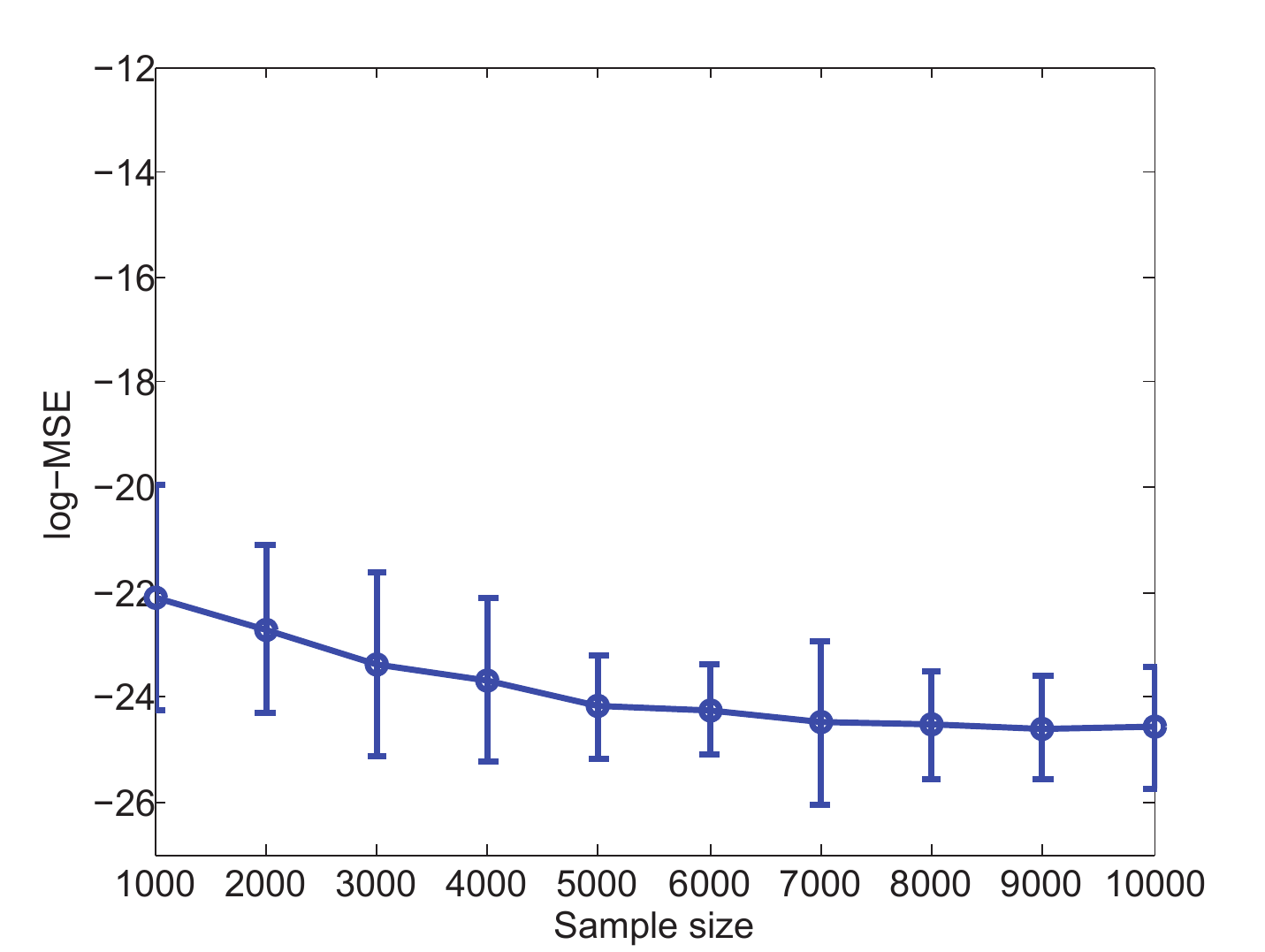} & 
\includegraphics[width=0.5\textwidth]{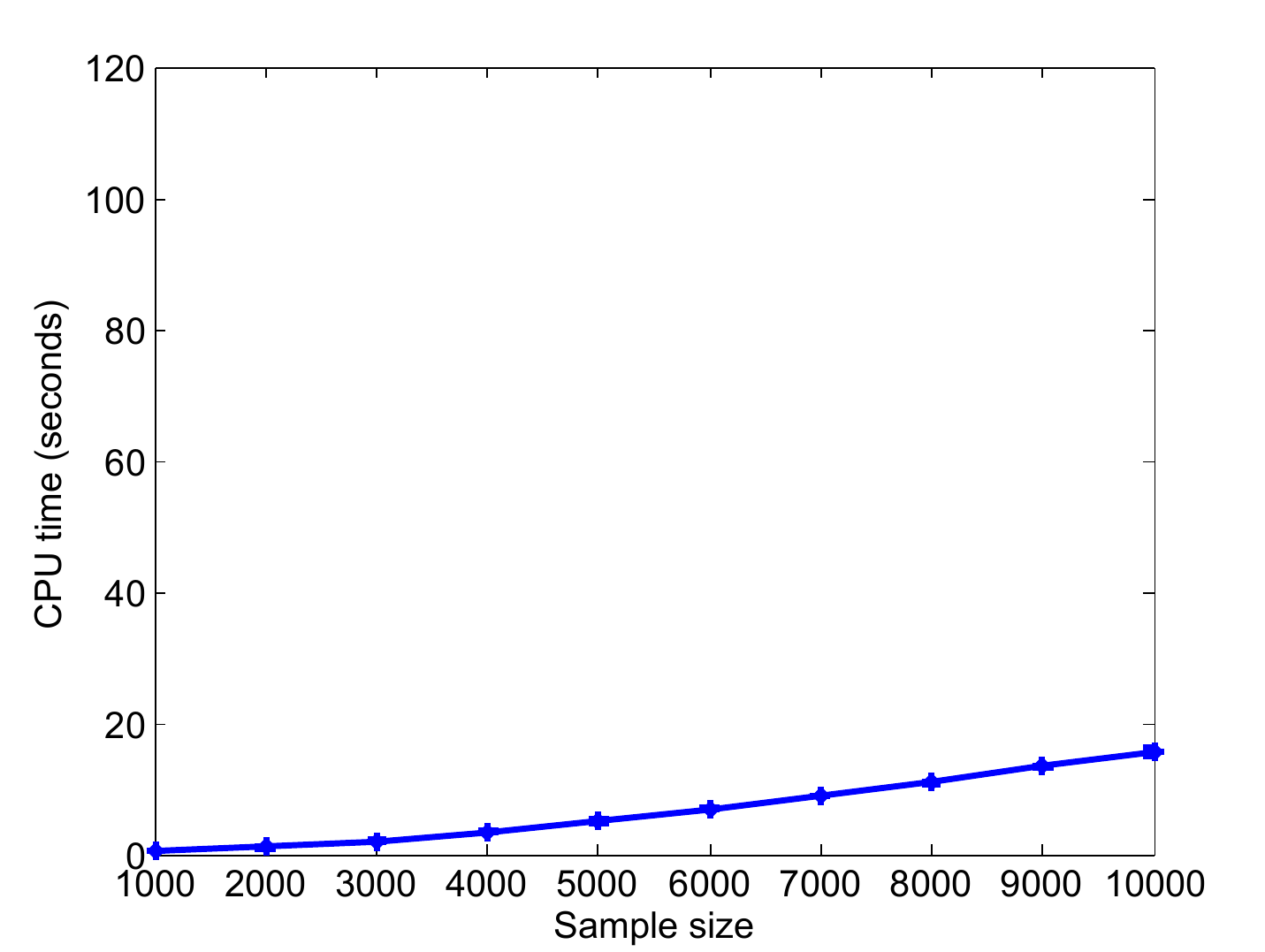} \\
\footnotesize{(a)} & \footnotesize{(b)} \\
\includegraphics[width=0.5\textwidth]{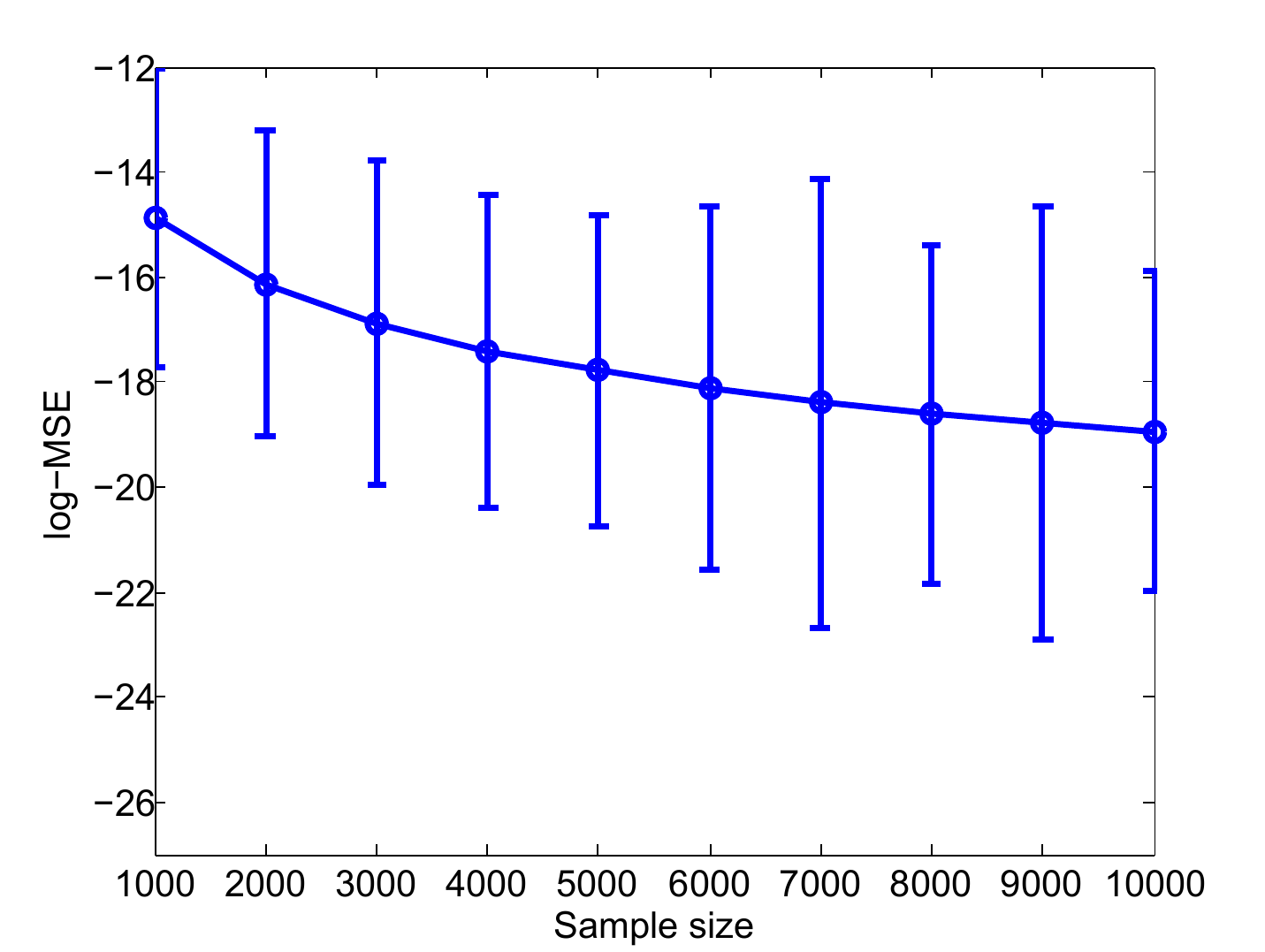} &
\includegraphics[width=0.5\textwidth]{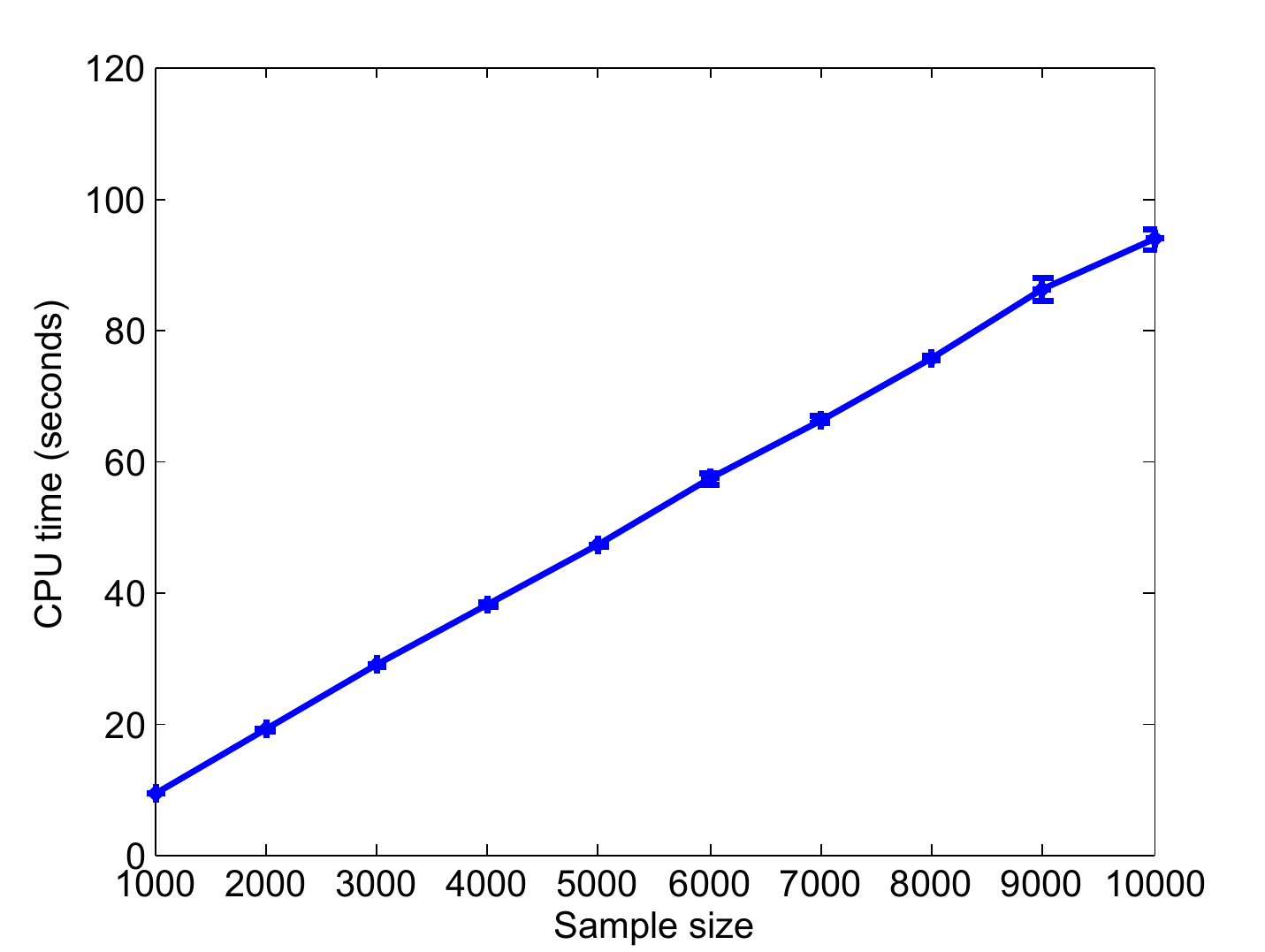} \\
 \footnotesize{(c)} & \footnotesize{(d)} \\
\includegraphics[width=0.5\textwidth]{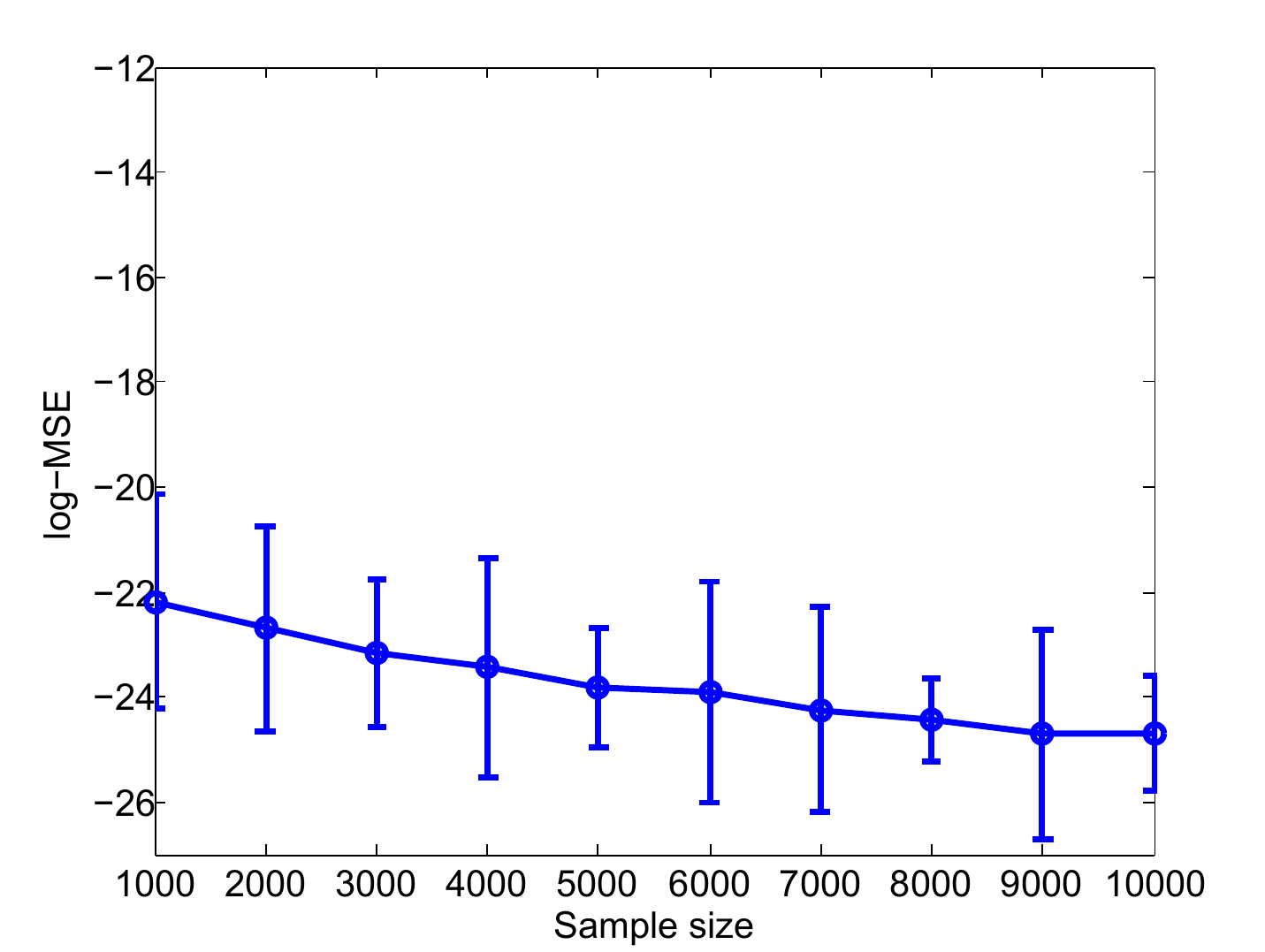} &
\includegraphics[width=0.5\textwidth]{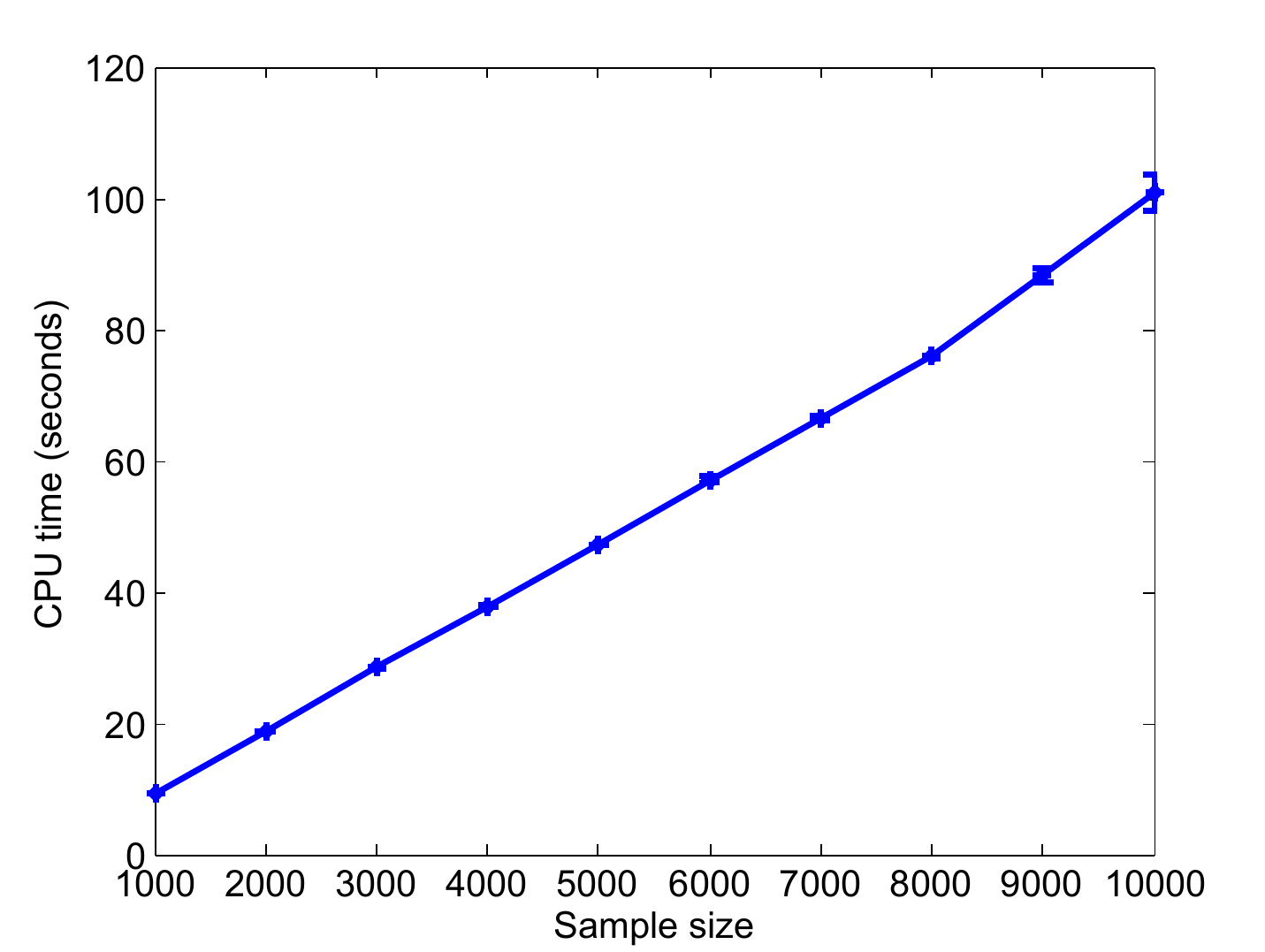} \\
 \footnotesize{(e)} & \footnotesize{(f)}
\end{tabular}}
\caption{\small Comparison of $\wh{f}_{\rm RWT}$ with the Ott and Kronmal (O\&K) estimator \cite{OttKro76} and with an exhaustive search for the best MSE. The plots show: (a) $\wh{f}_{\rm RWT}$ MSE; (b) $\wh{f}_{\rm RWT}$ times; (c) O\&K MSE; (d) O\&K times; (e) Exhaustive search MSE; (f) Exhaustive search times. All results are averaged over ten independent runs for each sample size (the error bars show the standard deviations). \label{fig:ottkronmal}}
\end{figure}

Additionally, we have performed comparisons between our proposed method and two alternatives: the Ott and Kronmal thresholding estimator \cite{OttKro76}; and an exhaustive search over all possible thresholds for the best MSE. As seen in Fig.~\ref{fig:ottkronmal}, our thresholding estimator provides close to the best possible MSE with far lower computational cost than the alternatives.

\subsection{High-dimensional simulations}
\label{sec:optim}

Although our algorithm has been implemented in MATLAB and therefore can be much further optimized for speed, we have devised several strategies for traversing the coefficient tree efficiently and circumventing computer-architecture related challenges for high-dimensional problems. We describe those strategies and demonstrate them using the same experimental set-up as above, but with much larger dimensionality.

\begin{itemize}

\item{\bf Direct computation of the coefficients near the leaves of the tree.} As discussed in \ref{ssec:sparsity}, the direct estimator (\ref{eq:weight_estimate_direct}) for $\wh{W}_u$ becomes computationally more efficient than the indirect estimator (\ref{eq:weight_estimate}) for  $k \ge d - \log (nd)$. Hence, near the bottom of the tree, instead of continuing to traverse the remaining levels based on the weights $\wh{W}_u$, we simply compute all coefficients $\wh{\theta}_u$ at the leaves of the corresponding subtree. This is always the most sensible course of action, given the fact that (\ref{eq:weight_estimate_direct}) requires the $\wh{\theta}_u$ anyway.

\item{\bf $N$-ary tree.} The number of levels to be traversed can be reduced by a factor of $\log N$ by considering an $N$-ary, rather than a binary, tree, at the cost of an increased number ($N$) of branches per level. We have found this trade-off to be worthwhile in many cases due to the possibility of vectorizing the computation of $\wh{W}_u$ for all the branches in each level and taking advantage of optimized routines for matrix algebra.

\item{\bf Open-node queue.} While our algorithm lends itself to recursive implementation, many computer/operating system architectures impose a hard limit on the recursion level due to stack-size restrictions (recursive function calls typically use stack memory). This becomes a problem when the dimension $d$ is high and the tree is accordingly very deep. We have circumvented the issue by implementing a queue system where the so-called open nodes (the tree branches awaiting processing) are sorted according to some criterion. This amounts to transferring the ``stack'' to user memory, where the only limit on the number of nodes is the free memory size. Possible criteria for sorting the nodes include depth-first, breadth-first and highest weight $\wh{W}_u$. We have used the latter criterion in our high-dimensional simulations.

\item {\bf Pruning high-frequency Walsh coefficients.} For a given Walsh function $\chr_s$, the Hamming weight  (i.e., number of ``on'' bits) of $s$ is a measure of frequency; higher-frequency coefficients have a higher proportion of ones in $s$. Because in many problems it is appropriate to assume that the signals of interest have low frequency, we have included the ability to impose a limit on the order of the Walsh coefficients by ignoring any branches with more than $m$ ``on" bits. The choice of $m$ depends on the problem context and on computational resources.

\item{\bf Weight-adaptive thresholding.} For some datasets, significant gains can be achieved by varying the thresholds $\lambda_{k,n}$ in a data-driven manner at each level of the tree; as an alternative to the preset schedules $\alpha_k$ in (\ref{eq:thresholds}), it is possible to take the weights $\wh{W}_u$ for each branch at level $k$, and then expand only the top $q$ branches with the highest $\wh{W}_u$. This is equivalent to making $\alpha_k$ not only level--dependent, but also dependent on the sequence of  weights at level $k$. The value of $q$ controls the trade-off between computation speed and accuracy.

\end{itemize}
The first three strategies are important modifications to a naive implementation of our algorithm, but in no way impact the MSE. The latter two techniques, however, provide an approximation to the estimator proposed and analyzed in this paper; for appropriate values of $m$ and $q$, they yield significant computational savings for a modest increase in MSE.

\begin{figure}[htbp]
\centerline{
\setlength{\tabcolsep}{-0.04cm}
\begin{tabular}{c} 
\includegraphics[width=0.6\textwidth]{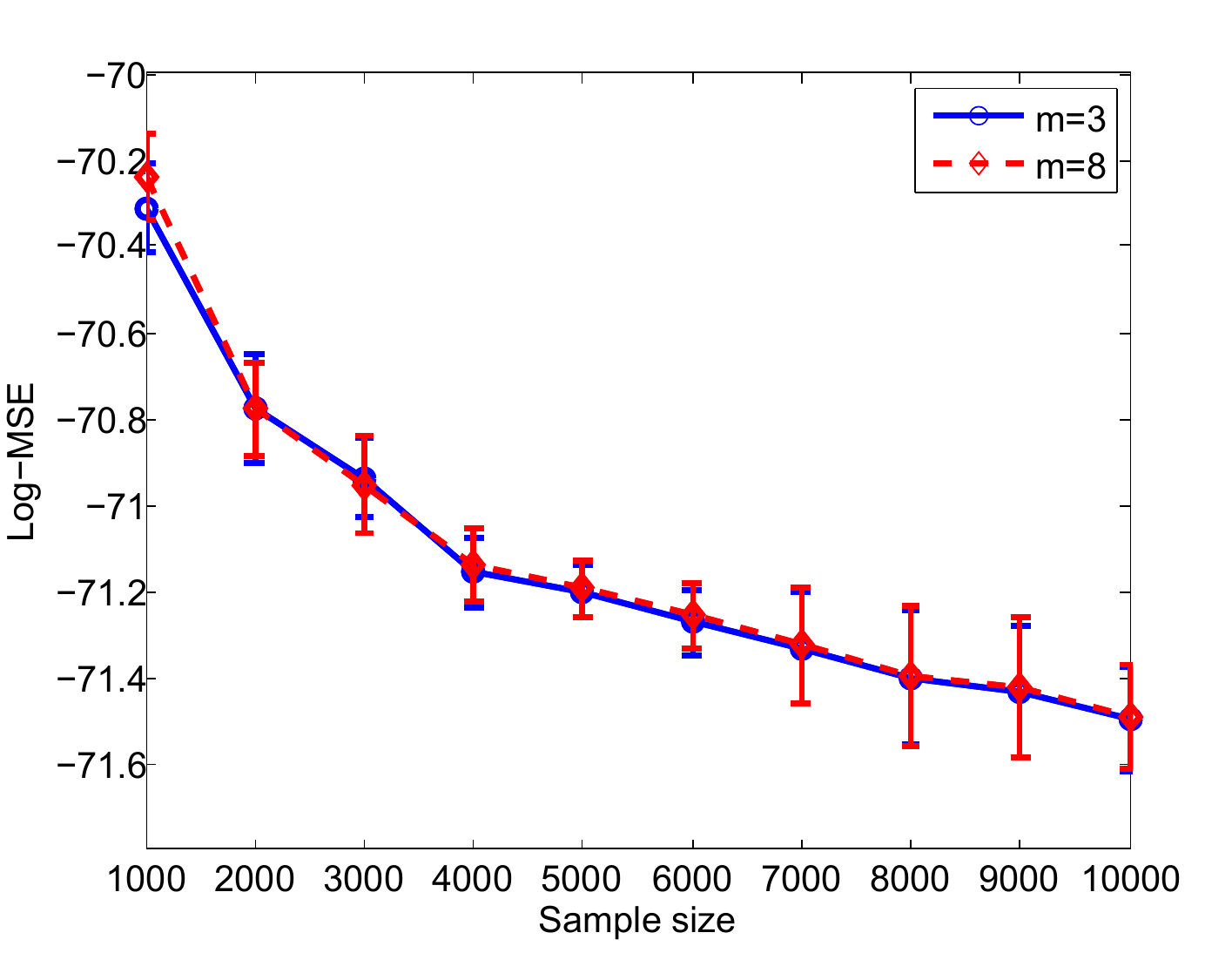} \\
\footnotesize{(a)} \\
\includegraphics[width=0.6\textwidth]{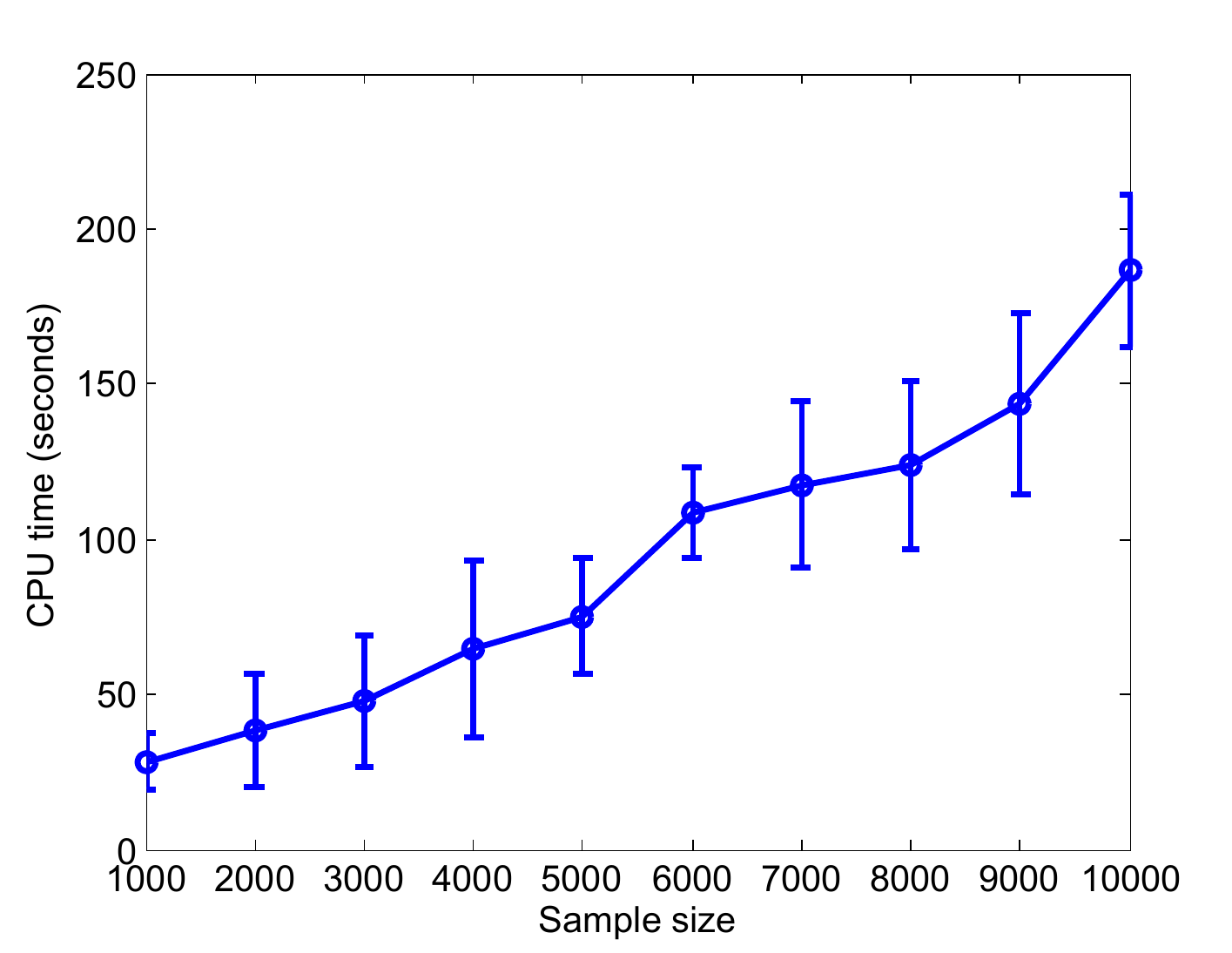} \\
\footnotesize{(b)} \\
\includegraphics[width=0.6\textwidth]{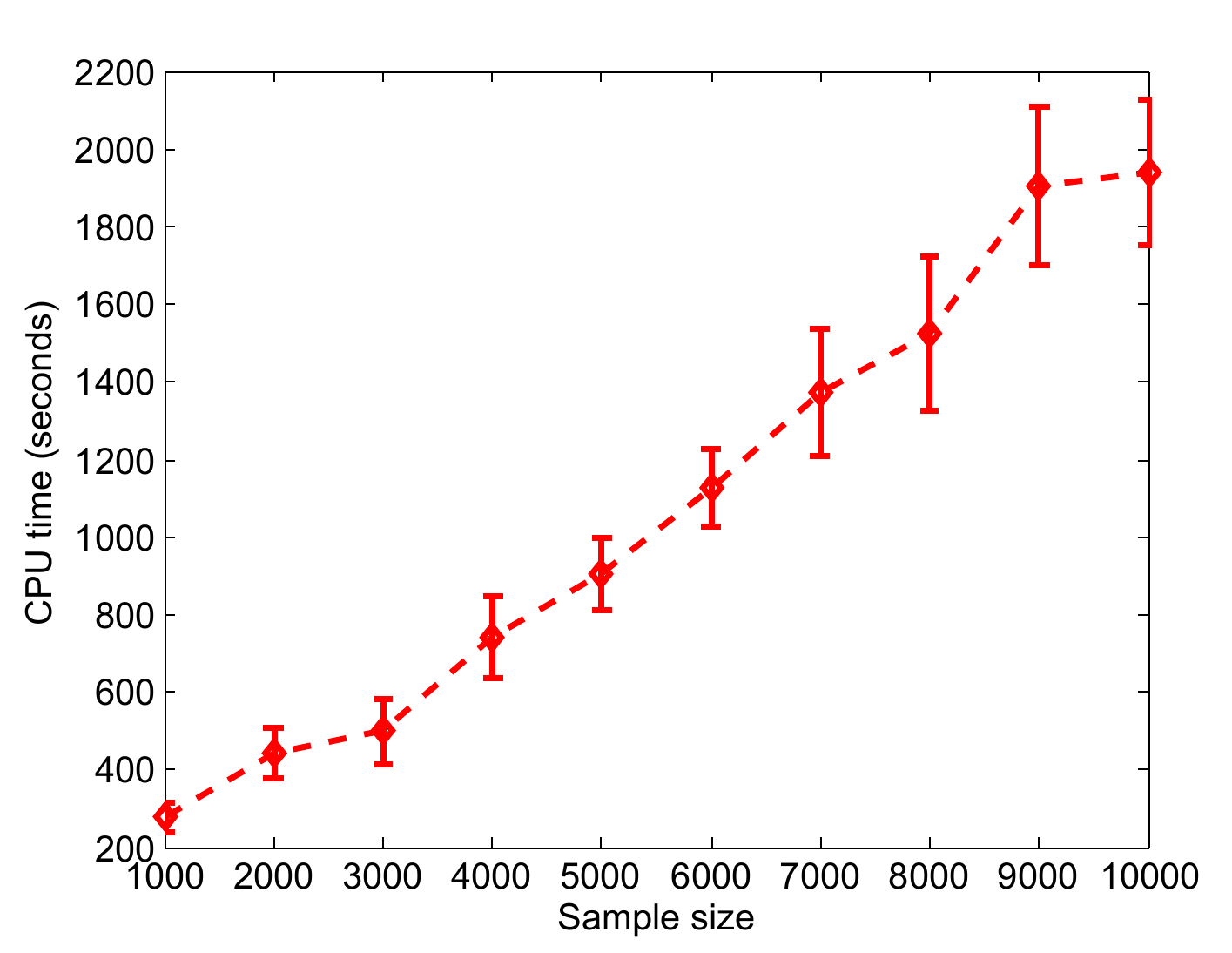} \\
\footnotesize{(c)}
\end{tabular}}
\caption{\small Performance for a large--dimensional problem ($d=50$): (a)~Log--MSE; (b)--(c)~running time (in seconds) for $m=3$ and $m=8$, respectively. All results are averaged over ten independent runs for each sample size (the error bars show the standard deviations). \label{fig:performance_large}}
\end{figure}

In Figure \ref{fig:performance_large}, we present plots of the MSE and computation time for simulated data with $d=50$ and multiple sample sizes ($n$), using the aforementioned optimizations. As before, the data were generated from a Bernoulli mixture density, similar to the one used for Figure \ref{fig:performance} but using ten $50$-dimensional mixture components, where each component has $47$ covariates with a {\rm Bernoulli}($1/2$) distribution and three covariates with a {\rm Bernoulli}($0.9$) distribution. The results are averaged over ten independent runs. We have limited the number of ``on'' bits in the Walsh binary strings to three and eight (\emph{i.e.}, $m=3$ and $m=8$ respectively), expanded only the $16$ subtrees with highest $\wh{W}_u$ at each level (\emph{i.e.}, $q=16$) and used an $N$-ary tree with $N=256$. The subtrees in the open-node queue were sorted by decreasing $\wh{W}_u$. Even in this high-dimensional regime, we achieve steadily decreasing MSE as a function of $n$, as well as approximately linear scaling in computation time. It is also apparent that setting $m=3$ achieves essentially the same MSE but with an order-of-magnitude reduction in the computational effort.

\section{Summary and conclusion}
\label{sec:conclusion}

We have presented a computationally efficient adaptive procedure for estimating a multivariate binary density in the ``big $d$, small $n$'' regime, which essentially forces a ``nonparametric'' approach. Many problems of current practical interest that involve multivariate binary data seem to pertain to populations with certain ``constellation'' effects among the $d$ covariates. We have formalized this observation by focusing on a class of densities whose Walsh representations exhibit a certain power-law behavior. For moderate sample sizes, our estimator attains nearly minimax rates of MSE convergence over this class and runs in polynomial time with high probability. Moreover, the complexity improves for sparser densities. We have also reported the results of simulations, which show that our implemented estimator behaves in accordance with the theory even in the small-sample regime.  In the future, we plan to test our method on real high-dimensional data sets. Another promising future direction is to investigate the relationship between various smoothness classes of densities on the binary hypercube defined in terms of their Fourier--Walsh representations and probability densities of binary Markov random fields \cite{Lauritzen}. Such a density will have the form
\begin{align*}
	f(x) &= \frac{1}{Z}\prod^d_{i=1}e^{-h_i(x^{(i)},x^{N_i})},
\end{align*}
where for each $i \in \{1,\ldots,d\}$ we have a neighborhood $N_i \subseteq \{1,\ldots,d\} \backslash \{i\}$, the corresponding ``local energy'' function $h_i(\cdot)$ depends only on $x^{(i)}$ and on $x^{N_i} = (x^{(j)} : j \in N_i)$, and $Z$ is the normalization constant known as the {\em partition function}. It is reasonable to assume that if most of the neighborhoods $N_i$ are small, then most of the Fourier--Walsh coefficients of $f$ will be small as well. Assuming specific bounds on the decay of the Fourier--Walsh coefficients of $f$ amounts to assuming something about the decay of correlations in the Markov random field governed by $f$. It remains to be seen whether sparsity classes of the type investigated in this paper can serve as a good model of binary Markov random fields with polynomial decay of correlations, or whether one would need to introduce a binary analog of something like (weak) Besov bodies \cite{BerLof76} in order to account for localization both in space (small number of $i$'s with large neighborhoods) and in frequency (small number of large Fourier--Walsh coefficients). In either case, it would be worthwhile to investigate the use of recursive thresholding estimators of the type introduced in this paper to estimate Markov graphical models on the binary hypercube.
 
\appendix
\renewcommand{\theequation}{\thesection.\arabic{equation}}
\setcounter{equation}{0}

\section{Auxiliary proofs}
\label{app:proofs}

\subsection{Proof of Lemma~\ref{lm:fu}}
\label{app:fu}

Using the appropriate definitions, as well as the factorization property (\ref{eq:walsh_factorize}) of the Walsh functions, we have
\begin{align*}
 \E_f \left\{ \chr_u(\pi_k(X)) I_{\{\sigma_k(X) = y\}}\right\} &= \sum_{x \in \bin^d} f(x) \chr_u(\pi_k(x))I_{\{\sigma_k(x)=y\}} \\
&= \sum_{z \in \bin^k} f(zy) \chr_u(z) \\
&= \sum_{z \in \bin^k} \Bigg(\sum_{(v,w) \in \bin^k \times \bin^{d-k}} \theta_{vw}\chr_{vw}(zy)\Bigg) \chr_u(z) \\
&= \sum_{z \in \bin^k} 
\Bigg(\sum_{(v,w) \in \bin^k \times \bin^{d-k}} \theta_{vw}\chr_v(z)\chr_w(y)\Bigg)\chr_u(z)  \\
&= \sum_{(v,w) \in \bin^k \times \bin^{d-k}} \theta_{vw}\chr_w(y) \ave{\chr_v,\chr_u} \\
&= \sum_{w \in \bin^{d-k}} \theta_{uw} \chr_w(y) \\
&\equiv f_u(y).
\end{align*}
Similarly,
\begin{align*}
W_u &= \sum_{y \in \bin^{d-k}} f^2_u(y) \\
&= \sum_{y \in \bin^{d-k}}  \E \left\{ \chr_u(\pi_k(X)) I_{\{\sigma_k(X) = y\}}\right\}  f_u(y) \\
&= \sum_{x \in \bin^d} \sum_{y \in \bin^{d-k}} f(x)\chr_u(\pi_k(x))f_u(y)I_{\{\sigma_j(x) = y\}} \\
&= \sum_{x \in \bin^d} f(x) \chr_u(\pi_k(x)) f_u(\sigma_k(x)) \\
&\equiv \E_f\left\{ \chr_u(\pi_k(X)) f_u(\sigma_k(X)) \right\},
\end{align*}
and the lemma is proved.

\subsection{Proof of Lemma~\ref{lm:wu}}
\label{app:wu}

We begin by showing that, for any $1 \le k < d$ and any $u \in \bin^k$,
\begin{equation}
\wh{W}_u = \wh{W}_{u0} + \wh{W}_{u1}.
\label{eq:W_additivity}
\end{equation}
From the factorization properties of the Walsh functions, we have
\begin{align*}
\wh{W}_{u0} &= \frac{1}{n^2}\sum^n_{i=1}\sum^n_{j=1} \chr_{u0}(\pi_{k+1}(X_i))\chr_{u0}(\pi_{k+1})X_j)) I_{\{\sigma_{k+1}(X_i) = \sigma_{k+1}(X_j)\}} \\
&= \frac{1}{2n^2}\sum^n_{i=1}\sum^n_{j=1}\chr_u(\pi_k(X_i))\chr_u(\pi_k(X_j))I_{\{\sigma_{k+1}(X_i) = \sigma_{k+1}(X_j)\}}
\end{align*}
and
\begin{align*}
\wh{W}_{u1} &= \frac{1}{n^2}\sum^n_{i=1}\sum^n_{j=1} \chr_{u1}(\pi_{k+1}(X_i))\chr_{u1}(\pi_{k+1})X_j)) I_{\{\sigma_{k+1}(X_i) = \sigma_{k+1}(X_j)\}} \\
&= \frac{1}{2n^2}\sum^n_{i=1}\sum^n_{j=1}\chr_u(\pi_k(X_i))\chr_u(\pi_k(X_j))(-1)^{X^{(k+1)}_i} \\
&\qquad \qquad \qquad \qquad \times (-1)^{X^{(k+1)}_l} I_{\{\sigma_{k+1}(X_i) = \sigma_{k+1}(X_j)\}}.
\end{align*}
Adding these two expressions and using the fact that
$$
1 + (-1)^{X^{(k+1)}_i}(-1)^{X^{(k+1)}_l} = 2I_{\{X^{(k+1)}_i = X^{(k+1)}_j\}},
$$
we get
\begin{align*}
\wh{W}_{u0} + \wh{W}_{u1} &= \frac{1}{n^2} \sum^n_{i=1}\sum^n_{j=1} \chr_u(\pi_k(X_i))\chr_u(\pi_k(X_j)) \\
& \qquad \qquad \qquad \qquad \times I_{\{\sigma_{k+1}(X_i) = \sigma_{k+1}(X_j), X^{(k+1)}_i = X^{(k+1)}_j\}} \\
&= \frac{1}{n^2} \sum^n_{i=1}\sum^n_{j=1} \chr_u(\pi_k(X_i))\chr_u(\pi_k(X_j)) I_{\{\sigma_k(X_i) = \sigma_k(X_j)\}} \\
&\equiv \wh{W}_u.
\end{align*}
This proves (\ref{eq:W_additivity}). By induction, we have
$$
\wh{W}_u = \sum_{s \in \cL(u)} \wh{W}_s,
$$
where $\cL(u)$ denotes the set of all leaves descending from $u$. Since $\cL(u) = \{uv \in \bin^d : v \in \bin^{d-k}\}$ and since
$$
\widehat{W}_s = \frac{1}{n^2}\sum^n_{i=1}\sum^n_{j=1}\chr_s(X_i)\chr_s(X_j) \equiv \wh{\theta}^2_s,
$$
the lemma is proved.

\subsection{Proof of Lemma~\ref{lm:separated}}

We first consider Item 1. By construction, for any $f \in \cF(k,a)$ we have
\begin{align*}
	f(x) &\le \frac{1}{M} + \sum^M_{j=2} |\theta_{s_j}(f)| \cdot \| \chr_{s_j} \|_\infty = \frac{1}{M} + \frac{ka}{\sqrt{M}}
\end{align*}
and
\begin{align*}
	f(x) &\ge \frac{1}{M}- \sum^M_{j=2} |\theta_{s_j}(f)| \cdot \| \chr_{s_j} \|_\infty = \frac{1}{M} - \frac{ka}{\sqrt{M}}
\end{align*}
Thus, if $k$ and $a$ satisfy the first condition in \eqref{eq:ka_conditions}, then all $f \in \cF(k,a)$ will be bounded between $1/2M$ and $3/2M$. To see that the second condition in \eqref{eq:ka_conditions} implies $\cF(k,a)\subset \cF_d(p)$, observe that in that case the Fourier--Walsh coefficients of $f$ ordered according to decreasing magnitude are
\begin{align*}
	|\theta_{(m)}(f)| &= \begin{cases}
	\frac{1}{\sqrt{M}}, & m=1 \\
	a, & m = 2,\ldots,k+1 \\
	0, & m = k+2,\ldots,M
\end{cases}
\end{align*}
from which it follows that $f \in \cF_d(p)$. Finally, any $f \in \cF(k,a)$ is also a probability density because it is nonnegative and because
$$
\sum_{x \in \bin^d} f(x) = \ave{f,1} = \sqrt{M}\ave{f,\chr_0} = \sqrt{M}\theta_0(f) = 1.
$$
This shows that $\cF(k,a) \subset \cF^{+,1}_d(p)$, as claimed.

We now move on to Item 2. To prove that the bound \eqref{eq:KL_bound} holds for any two densities $f,f' \in \cF(k,a)$, we use the fact that the Kullback--Leibler divergence $D(f \| f')$ is bounded above by the chi-square distance:
\begin{align}\label{eq:KL_vs_chi}
	D(f \| f') \le \chi^2(f,f') \deq \sum_{x \in \bin^d}\frac{\left|f(x)-f'(x)\right|^2}{f'(x)}.
\end{align}
For any $x \in \bin^d$ and $f,f' \in \cF(k,a)$, we have
\begin{align*}
	\left|f(x) - f'(x)\right|^2 &= \left|\sum^M_{j=1} \left(\theta_{s_j}(f) - \theta_{s_j}(f')\right)\chr_{s_j}(x) \right|^2 \le \frac{2ka^2}{M}.
\end{align*}
Using this together with the fact that $f,f' \ge 1/2M$ in \eqref{eq:KL_vs_chi}, we arrive at \eqref{eq:KL_bound}.

For Item 3, we need the following well-known combinatorial result from coding theory, the so-called Varshamov--Gilbert bound (see, e.g., Lemma~4.7 in \cite{Massart}): For any $m \in \N$, there exists a subset ${\cal K}_m$ of $\bin^m$ with the following properties:
\begin{itemize}
	\item for any $u,v \in {\cal K}_m$ with $u \neq v$,
	$$
	\sum^m_{i=1}I_{\{u^{(i)} \neq v^{(i)}\}} \ge m/4
	$$
	\item $\log |{\cal K}_m| \ge m/8$
\end{itemize}
There is a one-to-one correspondence between $\cF(M-1,a)$ and the binary hypercube $\bin^{M-1}$, under which each $f \in \cF(M-1,a)$ is mapped to $u_f \in \bin^{M-1}$ with
\begin{align*}
	u^{(j)}_f = I_{\{\theta_{s_{j+1}}(f) = a\}}, \qquad j = 1,\ldots,M-1.
\end{align*}
Given the Varshamov--Gilbert set ${\cal K}_{M-1} \subset \bin^{M-1}$, let
$$
\tilde{\cF}(M-1,a) \deq \Big\{ f \in \cF(M-1,a) : u_f \in {\cal K}_{M-1} \Big\}.
$$
Then for any two $f,f' \in \tilde{\cF}(M-1,a)$, we have
\begin{align}
	\| f - f' \|^2_{L^2(\mu_d)} &= 4a^2 \sum^{M-1}_{j=1} I_{\{u^{(j)}_f \neq u^{(j)}_{f'}\}} \ge (M-1)a^2,
\end{align}
and $\log |\cF(M-1,a)| = \log |{\cal K}_{M-1}| \ge (M-1)/8$. This takes care of Item 3.

Finally, we consider Item 4. To that end, we will need a refinement of the Varshamov--Gilbert bound due to Reynaud-Bouret \cite{ReynaudBouret}, which for our purposes can be stated as follows (see Lemma~4.10 in \cite{Massart}): For any $m,k \in \N$ with $m \ge k$, let $\bin^m_k$ denote the subset of $\bin^m$ that consists of all $u \in \bin^m$ with $\sum^m_{i=1} I_{\{u^{(i)}=1\}} = k$. If $m \ge 4k$, then there exists a set ${\cal K}_{m,k} \subset \bin^m_k$ with the following properties:
\begin{itemize}
	\item for any $u,v \in {\cal K}_{m,k}$ with $u \neq v$,
	$$
	\sum^m_{i=1}I_{\{u^{(i)}\neq v^{(i)}\}} \ge k/2
	$$
	\item $\log |{\cal K}_{m,k}| \ge 0.233 k \log (m/k)$
\end{itemize}
Now assume that $M-1 \ge 4k$ and consider the corresponding set ${\cal K}_{M-1,k}$. To each $u \in {\cal K}_{M-1,k}$ we can associate $2^k$ elements of $\cF(M-1,k)$, say $\{f_{u,v} :  v \in \bin^k \}$, such that $u$ determines the locations of the $k$ nonzero Walsh--Fourier coefficients taking values in $\{-a,a\}$, while the choice of $v \in \bin^k$ determines the signs of these coefficients. Thus, let
\begin{align*}
	\tilde{\cF}(k,a) \deq \Big\{ f_{u,v} : u \in {\cal K}_{m,k}, v \in \bin^k \Big\}.
\end{align*}
Consider now any two $f^{u,v},f_{u',v'} \in \tilde{\cF}(k,a)$, such that $(u,v) \neq (u',v')$. For any $j \in {1,\ldots,M-1}$ such that $u^{(j)} \neq u'^{(j)}$, either $u^{(j)}=0$ or $u'^{(j)}=0$. Suppose the latter. Then $|\theta_{s_{j+1}}(f_{u,v})| = a$, while $|\theta_{s_{j+1}}(f_{u',v'})| = 0$. Consequently,
\begin{align*}
	\left\| f_{u,v} - f_{u',v'} \right\|^2_{L^2(\mu_d)} &= a^2 \sum^{M-1}_{j=1} I_{\{u^{(j)}\neq u'^{(j)}\}} \ge ka^2
\end{align*}
and
\begin{align*}
	\log |\tilde{\cF}(k,a)| &= k \log 2 + \log |{\cal K}_{M-1,k}| \\
	&\ge k\log 2 + 0.233k \log \frac{M-1}{k} \\
	&\ge 0.233 k \left( \log \frac{M-1}{k} + 1\right).
\end{align*}
This proves Item 4.

\section{Empirical process representation}
\label{app:epr}
In this appendix, we show that for each $u \in \bin^{k}$, the
  norm $\| f_u - \wh{f}_u \|_{L^2(\mu_{d-k})}$ can be expressed as a
  supremum of an empirical process over a suitable function class;
  this was a key element of the proof of Theorem~\ref{thm:comp}.

\setcounter{equation}{0}

First, we show that $\| f_u - \wh{f}_u \|_{L^2(\mu_{d-k})}$ can be expressed as an empirical process of the form (\ref{eq:emp}) indexed by a suitable function class. To this end, define
$$
\cF \deq \left\{ \sum_{v \in \bin^{d-k}} \xi(v) \chr_{uv} : \xi \in L^2(\mu_{d-k}), \| \xi \|_{L^2(\mu_{d-k})} \le 1 \right\}.
$$
Then
\begin{equation}\label{eq:epr}
\| f_u - \wh{f}_u \|_{L^2(\mu_{d-k})} = \sup_{g \in \cF} \nu_n(g).
\end{equation}
Indeed, let $X_1,\ldots,X_n$ be $n$ i.i.d.~copies of $X \sim f$. Then
\begin{align*}
\nu_n(g) &= \frac{1}{n}\sum^n_{i=1} g(X_i) - \E g(X) \\
&= \frac{1}{n}\sum^n_{i=1} \sum_{v \in \bin^{d-k}} \xi(v) \chr_{uv}(X_i) - \sum_{v \in \bin^{d-k}} \xi(v)\theta_{uv} \\
&= \sum_{v \in \bin^{d-k}} (\wh{\theta}_{uv} - \theta_{uv}) \xi(v) \\
&\le \| \xi \|_{L^2(\mu_{d-k})} \| f_u - \wh{f}_u \|_{L^2(\mu_{d-k})},
\end{align*}
where in the last line we used Cauchy--Schwarz. This proves (\ref{eq:epr}). Next, we determine the constants $L$, $v$ and $H$ that are needed to apply Talagrand's bound (\ref{eq:talagrand}). For any $\xi \in L^2(\mu_{d-k})$ with unit norm, we have
\begin{align*}
\sup_{x \in \bin^d} \left|\sum_{v \in \bin^{d-k}} \xi(v) \chr_{uv}(x) \right| &\le \| \xi \|_{L^2(\mu_{d-k})} \sup_{x \in \bin^d} \sqrt{\sum_{v \in \bin^{d-k}} \chr^2_{uv}(x)} \\
& \le \sqrt{2^{d-k}2^{-d}} \\
&= 2^{-k/2}.
\end{align*}
Hence, any $g \in \cF$ is bounded by $L \equiv 2^{-k/2}$. From this, we also get the bound $\var g \le v$ with $v = L^2 = 2^{-k}$. Finally, to bound the Rademacher average (\ref{eq:rademacher}), we note that $\cF$ is the unit ball in the RKHS induced by the kernel
$$
K_u(x,y) \deq \sum_{v \in \bin^{d-k}} \chr_{uv}(x)\chr_{uv}(y), \qquad \forall x,y \in \bin^d.
$$
Then standard arguments (see, e.g., Section~2.4.2 in \cite{Men03}) lead to the bound
$$
\E \left\{ \sup_{g \in \cF} \sum^n_{i=1} \varepsilon_i g(U_i) \right\} \le \sqrt{\frac{n}{2^k}},
$$
which gives $H = 1/\sqrt{2^kn}$.

\bibliography{binary_hypercube.bbl}

\begin{thebibliography}{10}

\bibitem{AitAit76}
J.~Aitchison and C.~G.~G. Aitken.
\newblock Multivariate binary discrimination by the kernel method.
\newblock {\em Biometrika}, 63(3):413--420, 1976.

\bibitem{Bah61}
R.~R. Bahadur.
\newblock A representation of the joint distribution of $n$ dichotomous items.
\newblock In H.~Solomon, editor, {\em Studies in Item Analysis and Prediction},
  pages 169--176. Stanford Univ. Press, 1961.

\bibitem{BerLof76}
J.~Bergh and J.~L\"ofstr\"om.
\newblock {\em Interpolation {S}paces: {A}n {Introduction}}.
\newblock Springer-Verlag, 1976.

\bibitem{Can06}
E.~Cand\`es.
\newblock Modern statistical estimation via oracle inequalities.
\newblock {\em Acta Numerica}, 15:257--325, 2006.

\bibitem{CanTao06}
E.~J. Cand\`es and T.~Tao.
\newblock Near-optimal signal recovery from random projections: universal
  encoding strategies?
\newblock {\em IEEE Trans. Inform. Theory}, 52(12):5406--5425, December 2006.

\bibitem{Car07}
J.~M. Carro.
\newblock Estimating dynamic panel data discrete choice models with fixed
  effects.
\newblock {\em J. Econometrics}, 140:503--528, 2007.

\bibitem{CheKriLia89}
X.~R. Chen, P.~R. Krishnaiah, and W.~W. Liang.
\newblock Estimation of multivariate binary density using orthogonal functions.
\newblock {\em J. Multivariate Anal.}, 31:178--186, 1989.

\bibitem{CoverThomas}
T.~M. Cover and J.~A. Thomas.
\newblock {\em Elements of Information Theory}.
\newblock Wiley, New York, 2nd edition, 2006.

\bibitem{DFKO07}
I.~Dinur, E.~Friedgut, G.~Kindler, and R.~O'Donnell.
\newblock On the {F}ourier tails of bounded functions over the discrete cube.
\newblock {\em Israel J. Math.}, 160(389-412), 2007.

\bibitem{DJKP96}
D.~L. Donoho, I.~M. Johnstone, G.~Kerkyacharian, and D.~Picard.
\newblock Density estimation by wavelet thresholding.
\newblock {\em Ann. Statist.}, 24(2):508--539, 1996.

\bibitem{Efromovich}
S.~Efromovich.
\newblock {\em Nonparametric Curve Estimation}.
\newblock Springer, 1999.

\bibitem{dental}
M.~J. Garc\'{i}a-Zattera, A.~Jara, E.~Lesaffre, and D.~Declerck.
\newblock Conditional independence of multivariate binary data with an
  application in caries research.
\newblock {\em Computational Statistics and Data Analysis}, 51:3223--3234,
  2007.

\bibitem{GhaHel06}
Z.~Ghahramani and K.~Heller.
\newblock Bayesian sets.
\newblock In Y.~Weiss, B.~Sch\"{o}lkopf, and J.~Platt, editors, {\em Advances
  in Neural Information Processing Systems 18}, pages 435--442. MIT Press,
  Cambridge, MA, 2006.

\bibitem{gilbertMansour}
A.~C. Gilbert, S.~Guha, P.~Indyk, S.~Muthukrishnan, and M.~Strauss.
\newblock Near-optimal sparse {F}ourier representations via sampling.
\newblock In {\em Proceedings of the Thiry-Fourth Annual ACM Symposium on
  Theory of Computing}, 2002.

\bibitem{GilbertSPIE}
A.~C. Gilbert, S.~Muthukrishnan, and M.~J. Strauss.
\newblock Improved time bounds for near-optimal sparse {F}ourier representation
  via sampling.
\newblock In {\em Proc. {SPIE} Wavelets {XI}}, San Diego, CA, 2005.

\bibitem{GilbertStrauss}
A.~C. Gilbert and M.~J. Strauss.
\newblock Group testing in statistical signal recovery.
\newblock Preprint, 2006.

\bibitem{GolLev89}
O.~Goldreich and L.~Levin.
\newblock A hard-core predicate for all one-way functions.
\newblock In {\em Proc. 21st ACM Symp. on Theory of Computing}, pages 25--32,
  1989.

\bibitem{bacterialTaxonomy}
M.~Gyllenberg and T.~Koski.
\newblock Probabilistic models for bacterial taxonomy.
\newblock {\em International Statistical Review}, 69(2):249--276, August 2001.

\bibitem{HalKerPic98}
P.~Hall, G.~Kerkyacharian, and D.~Picard.
\newblock Block threshold rules for curve estimation using kernel and wavelet
  methods.
\newblock {\em Ann. Statist.}, 26(3):922--942, 1998.

\bibitem{HPKP97}
P.~Hall, S.~Penev, G.~Kerkyacharian, and D.~Picard.
\newblock Numerical performance of block thresholded wavelet estimators.
\newblock {\em Statistics and Computing}, 7:115--124, 1997.

\bibitem{Joh94}
I.~M. Johnstone.
\newblock Minimax {B}ayes, asymptotic minimax and sparse wavelet priors.
\newblock In S.~S. Gupta and J.~O. Berger, editors, {\em Statistical Decision
  Theory and Related Topics V}, pages 303--326. Springer, 1994.

\bibitem{KusMan93}
E.~Kushilevitz and Y.~Mansour.
\newblock Learning decision trees using the {F}ourier spectrum.
\newblock {\em {SIAM} J. Comput.}, 22(6):1331--1348, 1993.

\bibitem{Lauritzen}
S.~L. Lauritzen.
\newblock {\em Graphical Models}.
\newblock Clarendon Press, Oxford, 1996.

\bibitem{LiaKri85}
W.-Q. Liang and P.~R. Krishnaiah.
\newblock Nonparametric iterative estimation of multivariate binary density.
\newblock {\em J. Multivariate Anal.}, 16:162--172, 1985.

\bibitem{Man94}
Y.~Mansour.
\newblock Learning {B}oolean functions via the {F}ourier transform.
\newblock In V.~P. Roychodhury, K.-Y. Siu, and A.~Orlitsky, editors, {\em
  Theoretical {A}dvances in {N}eural {C}omputation and {L}earning}, pages
  391--424. Kluwer, 1994.

\bibitem{Massart}
P.~Massart.
\newblock {\em Concentration Inequalities and Model Selection}.
\newblock Springer, 2007.

\bibitem{Men03}
S.~Mendelson.
\newblock A few notes on statistical learning theory.
\newblock In S.~Mendelson and A.~J. Smola, editors, {\em Advanced Lectures in
  Machine Learning}, volume 2600 of {\em Lecture Notes in Computer Science}.
  Springer, 2003.

\bibitem{OttKro76}
J.~Ott and R.~A. Kronmal.
\newblock Some classification procedures for multivariate binary data using
  orthogonal functions.
\newblock {\em J. Amer. Stat. Assoc.}, 71(354):391--399, June 1976.

\bibitem{ReynaudBouret}
P.~Reynaud-Bouret.
\newblock Adaptive estimation of the intensity of inhomogeneous {P}oisson
  processes via concentration inequalities.
\newblock {\em Probab. Th. Rel. Fields}, 126:103--153, 2003.

\bibitem{Ros72}
H.~P. Rosenthal.
\newblock On the span in $l_p$ of sequences of independent random variables.
\newblock {\em Israel J. Math.}, 8:273--303, 1972.

\bibitem{ShmZha02}
I.~Shmulevich and W.~Zhang.
\newblock Binary analysis and optimization-based normalization of gene
  expression data.
\newblock {\em Bioinformatics}, 18(4):555--565, 2002.

\bibitem{silvaHypergraph}
J.~Silva and R.~Willett.
\newblock Hypergraph-based detection of anomalous high-dimensional
  co-occurrences.
\newblock {\em IEEE Trans. Pattern Anal. Mach. Intel.}, 31(3):563--569, 2009.

\bibitem{Sim95}
J.~S. Simonoff.
\newblock Smoothing categorical data.
\newblock {\em J. Statist. Planning and Inference}, 47:41--60, 1995.

\bibitem{Tal94a}
M.~Talagrand.
\newblock On {R}usso's approximate zero-one law.
\newblock {\em Ann. Probab.}, 22(3):1576--1587, 1994.

\bibitem{Tal94}
M.~Talagrand.
\newblock Sharper bounds for {G}aussian and empirical processes.
\newblock {\em Ann. Probab.}, 22:28--76, 1994.

\bibitem{TaoVu06}
T.~Tao and V.~H. Vu.
\newblock {\em Additive Combinatorics}.
\newblock Cambridge Univ. Press, 2006.

\bibitem{Vaart_Wellner}
A.~W. {van der Vaart} and J.~A. Wellner.
\newblock {\em Weak Convergence and Empirical Processes}.
\newblock Springer, 1996.

\bibitem{communityDNA}
J.~D. Wilbur, J.~K. Ghosh, C.~H. Nakatsu, S.~M. Brouder, and R.~W. Doerge.
\newblock Variable selection in high-dimensional multivariate binary data with
  applications to the analysis of microbial community {DNA} fingerprints.
\newblock {\em Biometrics}, 58:378--386, June 2002.

\bibitem{YangBarronTechRep}
Y.~Yang and A.~Barron.
\newblock Information-theoretic determination of minimax rates of convergence.
\newblock Technical Report~28, Department of Statistics, Iowa State University,
  1997.

\bibitem{YangBarron}
Y.~Yang and A.~Barron.
\newblock Information-theoretic determination of minimax rates of convergence.
\newblock {\em Ann. Statist.}, 27(5):1564--1599, 1999.

\end{thebibliography}

\end{document}